\documentclass[11pt,a4paper]{article}
\usepackage{microtype}
\usepackage{geometry}
\usepackage{amsthm,amsmath,amssymb,mathtools,bm}
\usepackage{enumitem}
\usepackage{xcolor,hyperref,fancyhdr,graphicx}

\hypersetup{
	colorlinks=true,
	linkcolor = blue,
	citecolor = red,
	pdftitle={Schatten estimates for additive Richardson iteration and GLT closure of the Karcher mean},
	pdfauthor={Samuele Ferri; Valerio Loi; Stefano Serra-Capizzano; Giacomo Tento}
}

\newtheorem{theorem}{Theorem}
\newtheorem{lemma}{Lemma}
\newtheorem{definition}{Definition}
\newtheorem{proposition}{Proposition}
\theoremstyle{remark}

\newcommand{\quoted}[1]{\textquotedblleft#1"}
\newcommand{\N}{\mathbb{N}}
\newcommand{\C}{\mathbb{C}}

\newcommand{\tr}{\operatorname{tr}}
\newcommand{\dd}{\,\mathrm{d}}

\newcommand{\Karcher}{\mathcal{G}}
\newcommand{\norm}[1]{\left\|#1\right\|}
\newcommand{\imagy}{\mathrm{i}}
\newcommand{\abs}[1]{\left|#1\right|}

\title{\bfseries Karcher means of structured matrices are asymptotically structured: a study via Schatten estimates for the additive Richardson iteration}
\author{%
	Samuele Ferri\textsuperscript{1},
	Valerio Loi\textsuperscript{1},
	Stefano Serra-Capizzano\textsuperscript{2,3,4},
	Giacomo Tento\textsuperscript{1}
	\\[0.7em]
	\small\textsuperscript{1}Department of Theoretical and Applied Sciences, University of Insubria, Varese, Italy\\
	\small\textsuperscript{2}Department of Science and High Technology, University of Insubria, Como, Italy\\
	\small\textsuperscript{3}CLAP research Center, University of Insubria, Como, Italy\\
	\small\textsuperscript{4}Department of Information Technology, Uppsala University, Uppsala, Sweden
}
\date{}

\begin{document}
	\maketitle
	
	\begin{abstract}
		In the present work, we first study the additive Richardson iteration for the Karcher mean of complex Hermitian positive definite (HPD) matrix sequences, by obtaining normalized Schatten estimates independent of the matrix dimension. These estimates are used to prove that Karcher means of structured matrices are asymptotically structured and this is proven for input HPD Toeplitz and generalized locally Toeplitz (GLT) matrix sequences. More precisely, we provide the a.c.s. approximation, which is the key for demonstrating the GLT closure under Schatten spectral bounds, and the GLT conclusion can be extended to matrix sequences with symbols that are positive definite almost everywhere.
		Ad hoc numerical experiments illustrate the spectral distribution consequences of these results and few related computational advantages. Concluding remarks and open problems are provided in the final section.
	\end{abstract}
	
	\section{Introduction}\label{sec:introduction}
	The Karcher mean gives the definition of the geometric mean of a finite number of Hermitian positive definite (HPD) matrices. When these matrices come from a discretization, their dimension grows as the mesh is refined, and a natural question is how the spectrum of their mean behaves in this limit.
	For structured matrix sequences and even for matrix sequences with hidden structure, the usual setting consists in considering the asymptotic spectral distribution in Weyl sense, which makes use of a measurable function, called a symbol \cite{BGr,GS17,Ty}. The question we consider here is whether taking the Karcher mean of a fixed number of generalized locally Toeplitz sequences \cite{GS17,GS18,GLT20-uni,GLT20}, that we abbreviate with GLT, preserves the GLT structure with GLT symbol equal to the pointwise Karcher mean of their symbols.
	For instance, given the Hermitian generating functions $f_1,\ldots,f_k$ which we assume positive definite almost everywhere and by considering the $s$-block $d$-level HPD Toeplitz matrices $T_{\bm n}(f_1),\ldots,T_{\bm n}(f_k)$, ${\bm n}=(n_1,\ldots,n_d)$ multi-index, we can prove that the sequence of the Karcher means $\{GM(T_{\bm n}(f_1),\ldots,T_{\bm n}(f_k))\}_{\bm n}$ is a $s$-block $d$-level HPD GLT matrix sequence whose GLT symbol is the Karcher mean of the symbols $GM(f_1,\ldots,f_k))$ i.e,
	\begin{displaymath}
		\{GM(T_{\bm n}(f_1),\ldots,T_{\bm n}(f_k))\}_{\bm n}\sim_{\mathrm{GLT}} GM(f_1,\ldots,f_k),
	\end{displaymath}
	which implies the eigenvalue (and singular value) distribution i.e.
	\begin{displaymath}
		\{GM(T_{\bm n}(f_1),\ldots,T_{\bm n}(f_k))\}_{\bm n}\sim_\lambda GM(f_1,\ldots,f_k).
	\end{displaymath}
	In other words, the matrix $GM(T_{\bm n}(f_1),\ldots,T_{\bm n}(f_k))$ is not Toeplitz in general, but
	\begin{displaymath}
		\{GM(T_{\bm n}(f_1),\ldots,T_{\bm n}(f_k))-T_{\bm n}(GM(f_1,\ldots,f_k))\}_{\bm n}
	\end{displaymath}
	is clustered to zero in both eigenvalue and singular value sense.
	
	In other words, $\{GM(T_{\bm n}(f_1),\ldots,T_{\bm n}(f_k))\}_{\bm n}$ is asymptotically close to the Toeplitz matrix sequence $\{T_{\bm n}(GM(f_1,\ldots,f_k))\}_{\bm n}$ and in this sense we give a meaning to the expression \quoted{Karcher means of structured matrices are asymptotically structured}. The above has also a practical numerical impact since the structured matrix $T_{\bm n}(GM(f_1,\ldots,f_k)))$ can be used as starting point in the Karcher iteration or even a simpler circulant or $\tau$ representative \cite{CN,NgBook,KoroSINUM,Sbit}, provided that the related circulant or $\tau$ matrix sequence is distributed as $GM(f_1,\ldots,f_k)$.

	\paragraph{Novelty of the results.}
	The discussion on Toeplitz geometric means applies to the much larger class of $s$-block $d$ level GLT $*$-algebras for any $s,d\ge 1$, also in the reduced version. See \cite{GS17,GS18,GLT20-uni,GLT20} for the standard GLT $*$-algebras, while the reduced version was started in \cite[pp. 398-399]{GLT03}, formally stated in \cite[Section 3.1.4]{GLT06} and was finally treated exhaustively in \cite{reduced}. Thus the applicability of our findings is broad since GLT matrix sequences include as special cases any approximation by local methods of integro-differential problems in multidimensional domain, where the local methods include Finite Differences, Finite Volumes, Finite Elements of any order, Discontinuous Galerkin, Isogeometric analysis with full or intermediate regularity; see the chapters and sections regarding these various applications in \cite{GS17,GS18,GLT20-uni,GLT20}.

	\paragraph{Technical content.}
	With reference to the technical content, we emphasize that one of the difficulty lies in passing between two limits, one regarding the matrix size and the other concerning the iteration parameter. At a fixed matrix size, the Karcher mean can be obtained as the limit of an iteration. While straightforward GLT calculus \cite{GS17,GS18,GLT20-uni,GLT20} can identify the symbol of any iterate (see the introduction of \cite{GEO25} and \cite{DEG25,TH}), we need precise information on the Schatten norm asymptotics to conclude the closure of the limit.
	The normalized Schatten estimates in this paper provide this control. More precisely, under uniform spectral bounds on the input sequences and the initial matrix, the iterates remain in the prescribed interval and satisfy normalized Schatten estimates independent of the matrix dimension. The proposed proof controls the second derivative along the affine segment of the additive step.
	Using these estimates, we provide a.c.s. approximations \cite{Tilli98,acs} which represent the key to prove GLT closure under Schatten spectral bounds, and the GLT conclusion can be extended to matrix sequences with symbols that are positive definite almost everywhere, as conjectured in \cite[pp. 13-14]{GEO25}.

	\paragraph{Matrix means and their computation.}
	Ando, Li and Mathias \cite{ALM04} formulated an axiomatic list of properties for multivariate matrix geometric means and it is the basis of our analysis. The Riemannian approach developed by Moakher \cite{Moakher05} and Bhatia and Holbrook \cite{BH06} defines the mean as the minimizer of a sum of squared distances on the positive definite cone.
	Lawson and Lim \cite{LL11} proved monotonicity of this least-squares mean and later treated positive definite operators through the Karcher equation \cite{LL14}. Lim and P\'alfia \cite{LP12} obtained the Karcher mean as a limit of matrix power means; their subsequent work gives existence and uniqueness for probability measures with a finite first moment \cite{LP18}.
	
	Bini and Iannazzo \cite[Theorem~2, p.~1703]{BI13} studied the additive Richardson iteration used here and proved local convergence under appropriate conditions on the step size and initial matrix. Other methods use Riemannian optimization; for example, Iannazzo and Porcelli \cite{IP18} combine a Barzilai--Borwein step with a nonmonotone line search. These methods provide the computational background. The distribution argument here requires a uniform normalized error estimate, which a local convergence theorem does not supply.

	\paragraph{Spectral distributions and related GLT results.}
	Locally Toeplitz sequences were introduced by Tilli \cite{Tilli98} and they represent the starting point of the larger framework of the GLT theory, which extends this calculus to broad maximal $*$-algebras of matrix sequences; see the monographs \cite{GS17,GS18} and the block theory in \cite{GLT20-uni,GLT20}, which is the specific setting used in this paper. Approximating classes of sequences, usually abbreviated as a.c.s. \cite{acs}, allow approximation by errors of small norm together with errors of small relative rank and define, and provide the meaningful notion of convergence used in asymptotic spectral analysis of GLT sequences. The relation between this convergence in the space of matrix sequences and convergence in measure of symbols is developed in \cite{BG17,GLT20}.
	
	For geometric means of GLT sequences, in \cite[Section 3.1, Theorems~3--5]{GEO25} it is proved the two-input closure results for scalar and block multilevel sequences when at least one input symbol is invertible almost everywhere. Their discussion of more than two inputs \cite[Section 3.2]{GEO25} formulates the Karcher closure problem we analyse in this paper and identifies a.c.s. approximation by fixed iterates as the step needed to pass to the limit. In \cite[Theorems~3 and~5]{DEG25} the invertibility assumption in the two-input case for commuting symbols is removed; see also \cite{TH}.
	The finite-tuple problem with positive definite symbols considered here is therefore more general and requires substantially different tools from the two-input problem and from the case of symbols that are singular on a set of positive measure.

	\paragraph{The argument and its scope.}
	In Theorem~\ref{th:task} we state dimension-independent normalized Schatten estimates for the additive Richardson iteration under common bounds $aI\preceq A_{\bm n}^{(i)},X_{\bm n,0}\preceq bI$, with $0<a\le b<\infty$. The argument first keeps the iterates in this interval, then estimates the objective along the affine segment of the additive step. Theorem~\ref{th:task} supplies an error estimate that gives an a.c.s. approximation with the iteration count as approximation parameter. Under the same uniform bounds, with the GLT algebra axioms we identify the symbol of every fixed iterate and the convergence of these symbols and GLT closure then give the symbol of the Karcher mean.
	
	We further remove the uniform spectral bounds by placing the mean between two GLT approximations and provide the closure statement of Theorem~\ref{thm:glt-general}. 
	
	\paragraph{Organization of the work.}
	In Section~\ref{sec:preliminaries} we set the notation and we report the tools in axiomatic flavor as ACS1--ACS6, GLT1--GLT9, and K1--K11. The Schatten estimate and its detailed argument follow in Section~\ref{sec:schatten-main}. In Section~\ref{sec:glt-karcher} we treat the GLT closure in the bounded case, the squeezing lemma, and the general GLT closure theorem, before the experimental Section~\ref{sec:numerical-glt} and the conclusion Section~\ref{sec:end}.

	\section{Preliminaries}\label{sec:preliminaries}
	In the present section first we set the notation and then we list and discuss the various tools we use: distribution theory, approximating classes of sequences, GLT $*$-algebras, the Karcher mean and the axiomatic treatment of a matrix geometric mean.
	
	\subsection{Notation}
	We follow the notation of \cite[Sections 2.1.1--2.1.2, pp.~116--118]{GLT20}.
	Given the number of levels $d$, the block size $s$, and the number of inputs $k$, we set
	\[
	\bm n=(n_1,\ldots,n_d)\in\N^d,
	\qquad
	N(\bm n)=\prod_{j=1}^{d} n_j,
	\]
	with matrix size equal to $sN(\bm n)$. The notation $\bm n\to\infty$ means $\min_j n_j\to\infty$. The norm $\|\cdot\|_p$ is the Schatten norm of index $p$ and, in particular, $\|\cdot\|_2$ is the Frobenius norm. We write $\|\cdot\|$ for the spectral norm which is the the Schatten infinity norm. For a matrix $B$ of order $m$,
	\[
	\|B\|_p=\left(\sum_{j=1}^m\sigma_j(B)^p\right)^{1/p}
	\quad(1\le p<\infty),\qquad
	\|B\|=\max_j\sigma_j(B),
	\]
	where $\sigma_j(B)$, $j=1,\ldots,m$, are the singular values of $B$. We denote by $I$ and $O$ the identity and zero matrices of the appropriate size, while $\preceq$ is the L\"owner order.
	
	Throughout, $|D|$ denotes the uniform Lebesgue measure of a measurable set $D$. The number of inputs $k$ is fixed as the dimension grows.

	\subsection{Asymptotic spectral distribution}\label{sec:spectral-prelim}
	Write $m_{\bm n}=sN(\bm n)$ for the total matrix size. Let $D\subset\mathbb R^q$ be measurable with $0<|D|<\infty$, and let $f:D\to\C^{s\times s}$ be a measurable matrix valued function. Write $\sigma_j(B)$ and $\lambda_j(B)$ for the singular values and eigenvalues of a matrix $B$, possibly counted with multiplicity. We say that a sequence $\{A_{\bm n}\}_{\bm n}$, with matrices of size $m_{\bm n}$, has singular-value distribution $f$, and write $\{A_{\bm n}\}_{\bm n}\sim_\sigma f$, if
	\begin{equation}\label{eq:def-sv-distribution}
		\lim_{\bm n\to\infty}\frac1{m_{\bm n}}
		\sum_{j=1}^{m_{\bm n}}F(\sigma_j(A_{\bm n}))
		=\frac1{s|D|}\int_D\sum_{j=1}^sF(\sigma_j(f(z)))\,\mathrm dz
		\qquad(F\in C_c(\mathbb R)).
	\end{equation}
	Likewise, $\{A_{\bm n}\}_{\bm n}\sim_\lambda f$ means
	\begin{equation}\label{eq:def-eigenvalue-distribution}
		\lim_{\bm n\to\infty}\frac1{m_{\bm n}}
		\sum_{j=1}^{m_{\bm n}}F(\lambda_j(A_{\bm n}))
		=\frac1{s|D|}\int_D\sum_{j=1}^sF(\lambda_j(f(z)))\,\mathrm dz
		\qquad(F\in C_c(\C)).
	\end{equation}
	Here $C_c$ denotes the continuous compactly supported test functions; see \cite[Definition~2.11]{GLT20}. For Hermitian matrices and Hermitian $f$, it is enough to use test functions on $\mathbb R$.
	
	A sequence $\{Z_{\bm n}\}_{\bm n}$ is \emph{zero-distributed} if $\{Z_{\bm n}\}_{\bm n}\sim_\sigma0$. Equivalently, for every $\varepsilon>0$, the fraction of singular values exceeding $\varepsilon$ tends to zero. It can also be characterized by $Z_{\bm n}=R_{\bm n}+E_{\bm n}$ with $\operatorname{rank}(R_{\bm n})/m_{\bm n}\to0$ and $\|E_{\bm n}\|\to0$; see \cite[Theorem~2.18]{GLT20}. In particular, a low-rank part needs not have small or even controlled norm.
	
	Notice that $d$-level Toeplitz matrix sequences are characterized by $D=[-\pi,\pi]^d$ so $q=d$, while $d$-level GLT matrix sequences are characterized by $D=[0,1]^d\times [-\pi,\pi]^d$ so $q=2d$, because GLT matrix sequences include variation in space or in time or in both, as required in concrete applications when approximating spatially-variant integro-differential equations.

	\subsection{Approximating classes of sequences}\label{sec:acs-prelim}
	Approximating classes of sequences, abbreviated a.c.s., separate an approximation parameter from the matrix dimension. For each fixed approximation parameter, the dimension is allowed to grow. The error bounds are then sent to zero by increasing the approximation parameter.
	The theory permits both small-norm errors and errors confined to a vanishing fraction of the dimensions \cite[Section 2.5]{GLT20}.
	
	\begin{definition}[Approximating class of sequences]\label{def:acs}
		Let $A_{\bm n}$ and $B_{\bm n,r}$ have order $m_{\bm n}$. We write $\{B_{\bm n,r}\}_{\bm n}\xrightarrow{\mathrm{a.c.s.}} \{A_{\bm n}\}_{\bm n}$ as $r\to\infty$ if, for every $r\in\N$, there is $n_r$ such that, whenever $\min_j n_j\ge n_r$,
		\[
		A_{\bm n}=B_{\bm n,r}+R_{\bm n,r}+E_{\bm n,r},\qquad
		\operatorname{rank}(R_{\bm n,r})\le c(r)m_{\bm n},\qquad
		\|E_{\bm n,r}\|\le\omega(r),
		\]
		where $c(r)\to0$ and $\omega(r)\to0$ as $r\to\infty$.
		The functions $c,\omega$ and the threshold $n_r$ are independent of $\bm n$ \cite[Definition~2.31]{GLT20}.
	\end{definition}
	
	The squeezing argument below measures approximation errors by a bounded singular-value gauge function. Set $\psi(t)=t/(1+t)$ for $t\ge0$, and define
	\[
	p_\psi(B)=\frac1{\ell}\sum_{j=1}^{\ell}\psi(\sigma_j(B))
	\quad(B\in\C^{\ell\times\ell}),\qquad
	p_\psi^{\mathrm{acs}}(\{A_{\bm n}\})
	=\limsup_{\bm n\to\infty}p_\psi(A_{\bm n}),
	\]
	\begin{equation}\label{eq:acs-pseudometric-prelim}
		d_\psi(\{A_{\bm n}\},\{B_{\bm n}\})
		=\limsup_{\bm n\to\infty}p_\psi(A_{\bm n}-B_{\bm n}).
	\end{equation}
	This is a complete pseudometric on sequences of the prescribed dimensions.
	Its zero-distance pairs differ by a zero-distributed sequence, and Definition~\ref{def:acs} is equivalent to $d_\psi(\{A_{\bm n}\},\{B_{\bm n,r}\})\to0$
	\cite{BG17} and \cite[Theorem~2.32]{GLT20}.
	For symbols on $D$, the corresponding quantities are
	\[
	p_\psi^{\mathrm{meas}}(f)
	=\frac1{s|D|}\int_D\sum_{j=1}^s\psi(\sigma_j(f(z)))\,\mathrm dz,
	\qquad
	d_\psi^{\mathrm{meas}}(f,g)=p_\psi^{\mathrm{meas}}(f-g).
	\]
	Then $d_\psi^{\mathrm{meas}}(f_r,f)\to0$ if and only if $f_r\to f$ in measure, meaning $|\{z\in D:\|f_r(z)-f(z)\|>\varepsilon\}|\to0$ for every $\varepsilon>0$ \cite[Section 2.3.4]{GLT20}.
	
	For completeness, we list the standard properties ACS1--ACS6, keeping the numbering of \cite[Section 5, p.~179]{GLT20}.
	Recall that a sequence is called \emph{sparsely unbounded} if
	\[
	\lim_{M\to\infty}\limsup_{\bm n\to\infty}
	\frac{\#\{j:\sigma_j(A_{\bm n})>M\}}{m_{\bm n}}=0.
	\]
	Every sequence with a finite measurable singular-value symbol has this property \cite[Proposition~2.22]{GLT20}.
	\begin{enumerate}[label=\textbf{ACS\arabic*.},ref=ACS\arabic*,leftmargin=*]
		\item\label{acs:gauge} \textbf{Gauge identity.}
		If $\{A_{\bm n}\}_{\bm n}\sim_\sigma f$, then $p_\psi^{\mathrm{acs}}(\{A_{\bm n}\})=p_\psi^{\mathrm{meas}}(f)$ \cite[Theorem~2.33]{GLT20}.
		
		\item\label{acs:sv-transfer} \textbf{Transfer of singular-value distributions.}
		Suppose that $\{B_{\bm n,r}\}_{\bm n}\sim_\sigma f_r$ for every $r$, and
		\[
		\{B_{\bm n,r}\}_{\bm n}\xrightarrow{\mathrm{a.c.s.}}
		\{A_{\bm n}\}_{\bm n},\qquad f_r\to f\ \text{in measure}.
		\]
		Then $\{A_{\bm n}\}_{\bm n}\sim_\sigma f$. \cite[Theorem~2.34]{GLT20}.
		
		\item\label{acs:eigen-transfer} \textbf{Transfer of eigenvalue distributions.}
		The preceding statement holds with $\sim_\lambda$ in place of 	$\sim_\sigma$ when all the matrix sequences considered are Hermitian \cite[Theorem~2.35]{GLT20}.
		
		\item\label{acs:algebra} \textbf{Algebraic closure.}
		Suppose $B_{\bm n,r}\xrightarrow{\mathrm{a.c.s.}}A_{\bm n}$ and $C_{\bm n,r}\xrightarrow{\mathrm{a.c.s.}}D_{\bm n}$. Then adjoints and linear combinations 	converge in the same way. Also $B_{\bm n,r}C_{\bm n,r}\xrightarrow{\mathrm{a.c.s.}}A_{\bm n}D_{\bm n}$ if the two limit sequences are s.u. Right multiplication by a fixed sparsely unbounded sequence preserves a.c.s. convergence \cite[Theorem~2.37]{GLT20}.
		
		\item\label{acs:blocks} \textbf{Blocking property.}
		Entrywise a.c.s. convergence in a fixed-size block partition implies a.c.s. convergence of the global matrix sequences \cite[Theorem~2.38]{GLT20}.
		
		\item\label{acs:schatten} \textbf{Schatten criterion.}
		For any fixed $1\le p<\infty$, the condition
		\[
		\lim_{r\to\infty}\limsup_{\bm n\to\infty}
		\frac{\|A_{\bm n}-B_{\bm n,r}\|_p^p}{m_{\bm n}}=0
		\]
		implies a.c.s. convergence. For $p=\infty$, the sufficient condition is
		\[
		\lim_{r\to\infty}\limsup_{\bm n\to\infty}
		\|A_{\bm n}-B_{\bm n,r}\|=0.
		\]
		See \cite[Theorem~2.39 and Section 5, ACS6]{GLT20}.
	\end{enumerate}
	The Schatten criterion is the connection of our asymptotic distribution goal with the dimension-independent Richardson estimates.

	\subsection{Multilevel block GLT sequences}\label{sec:glt-prelim}
	In the theory of GLT sequences, we keep track of an associated canonical symbol that gives us information about both the spatial variable and the Fourier behavior of the discretized operator. The canonical symbols' domain is
	\[
	\Omega=[0,1]^d\times[-\pi,\pi]^d,\qquad
	z=(\bm{x},\bm\vartheta),\qquad |\Omega|=(2\pi)^d.
	\]
	A $d$-level $s$-block GLT sequence is composed by matrices of order $m_{\bm n}$ and a canonical measurable symbol $\kappa:\Omega\to\C^{s\times s}$; we write $\{A_{\bm n}\}_{\bm n}\sim_{\mathrm{GLT}}\kappa$.
	The classical constructive definition starts from finite sums of products of diagonal sampling sequences and block Toeplitz sequences, and closes the resulting matrix-symbol pairs under a.c.s. convergence in the matrices and convergence in measure in the symbols \cite[Definition~4.5 and Remark~4.21]{GLT20}.
	Unlike a singular-value symbol, which is invariant up to measure-preserving transformations, the canonical GLT symbol is unique up to equality in measure.
	
	For the basic generators, let
	\[
	\widehat f_j=(2\pi)^{-d}\int_{[-\pi,\pi]^d}
	f(\bm\vartheta)e^{-\mathrm i j\cdot\bm\vartheta}\,\mathrm d\bm\vartheta,
	\qquad
	T_{\bm n}(f)=[\widehat f_{i-j}]_{i,j=\bm1}^{\bm n},\qquad
	D_{\bm n}(a)=\operatorname{diag}_{j=\bm1,\ldots,\bm n}a(j/\bm n),
	\]
	where multi-indices have a lexicographical order and $j/\bm n$ is the pointwise division. The following properties retain the standard GLT1--GLT9 numbering of \cite[Section 5, pp.~179--180]{GLT20}, and they suffice to characterize GLT sequences for our purposes.
	\begin{enumerate}[label=\textbf{GLT\arabic*.},ref=GLT\arabic*,leftmargin=*]
		\item\label{glt:distribution} \textbf{Spectral distribution.}
		If $\{A_{\bm n}\}_{\bm n}\sim_{\mathrm{GLT}}\kappa$, then $\{A_{\bm n}\}_{\bm n}\sim_\sigma\kappa$. If every $A_{\bm n}$ is Hermitian, then $\kappa$ is Hermitian almost everywhere and $\{A_{\bm n}\}_{\bm n}\sim_\lambda\kappa$ \cite[Theorems~4.7 and~4.11, Proposition~4.10]{GLT20}.
		
		\item\label{glt:almost-hermitian} \textbf{Almost-Hermitian sequences.}
		Suppose $\{A_{\bm n}\}_{\bm n}\sim_{\mathrm{GLT}}\kappa$ and $A_{\bm n}=H_{\bm n}+E_{\bm n}$, where $H_{\bm n}$ is Hermitian and $m_{\bm n}^{-1/2}\|E_{\bm n}\|_2\to0$. Then $\{P_{\bm n}^*A_{\bm n}P_{\bm n}\}_{\bm n} \sim_\sigma\kappa$ and $\sim_\lambda\kappa$ whenever $P_{\bm n}\in\C^{m_{\bm n}\times\ell_{\bm n}}$, $P_{\bm n}^*P_{\bm n}=I_{\ell_{\bm n}}$, and $\ell_{\bm n}/m_{\bm n}\to 1$\cite[Theorem~4.12]{GLT20}. Taking $P_{\bm n}=I$ gives the uncompressed assertion. For the compressed sequence, the normalization in \eqref{eq:def-sv-distribution}--\eqref{eq:def-eigenvalue-distribution} uses its order $\ell_{\bm n}$.
		
		\item\label{glt:generators} \textbf{Basic sequences.}
		The symbol of $\{T_{\bm n}(f)\}_{\bm n}$ is $f(\bm\vartheta)$ for $f\in L^1([-\pi,\pi]^d;\C^{s\times s})$.
		The symbol of $\{D_{\bm n}(a)\}_{\bm n}$ is $a(\bm{x})$ when $a:[0,1]^d\to\C^{s\times s}$ is continuous almost everywhere.
		Finally, $\{Z_{\bm n}\}_{\bm n}\sim_{\mathrm{GLT}}O_s$ if and only if it is zero-distributed \cite[Theorems~3.11,~3.15 and~4.23]{GLT20}.
		
		\item\label{glt:algebra} \textbf{Algebra and inversion.}
		If $\{A_{\bm n}\}\sim_{\mathrm{GLT}}\kappa$ and $\{B_{\bm n}\}\sim_{\mathrm{GLT}}\xi$, then
		\[
		\{A_{\bm n}^*\}\sim_{\mathrm{GLT}}\kappa^*,\qquad
		\{\alpha A_{\bm n}+\beta B_{\bm n}\}\sim_{\mathrm{GLT}}
		\alpha\kappa+\beta\xi,\qquad
		\{A_{\bm n}B_{\bm n}\}\sim_{\mathrm{GLT}}\kappa\xi
		\]
		for $\alpha,\beta\in\C$. Product order is preserved.
		If $\kappa$ is invertible almost everywhere, then $\{A_{\bm n}^{\dagger}\}\sim_{\mathrm{GLT}}\kappa^{-1}$, where $\dagger$ is the Moore-Penrose inverse. When all finite matrices are invertible, this is the usual inverse \cite[Theorems~4.16--4.17]{GLT20}.
		
		\item\label{glt:functional-calculus} \textbf{Continuous functional calculus.}
		If $\{A_{\bm n}\}\sim_{\mathrm{GLT}}\kappa$ and every $A_{\bm n}$ is Hermitian, then
		$\{h(A_{\bm n})\}\sim_{\mathrm{GLT}}h(\kappa)$ for every continuous $h:\mathbb R\to\C$ \cite[Theorem~4.18]{GLT20}.
		A function continuous on a common compact spectral interval may be used after extending it continuously to $\mathbb R$. This observation will justify logarithms under uniform positive spectral bounds.
		
		\item\label{glt:block-assembly} \textbf{Block assembly.}
		For $1\le i,j\le q$, suppose $\{A_{\bm n,ij}\}_{\bm n}$ is $d$-level $s$-block GLT with symbol $\kappa_{ij}$. Let $A_{\bm n}=[A_{\bm n,ij}]_{i,j=1}^q$. Define the permutation $\Pi_{\bm n,q,r}$ by
		\[
		\Pi_{\bm n,q,r}(e_i\otimes e_j\otimes e_\alpha)
		=e_j\otimes e_i\otimes e_\alpha,
		\quad 1\le i\le q,\quad1\le j\le N(\bm n),\quad1\le\alpha\le r.
		\]
		Then $\{\Pi_{\bm n,q,r}A_{\bm n}\Pi_{\bm n,q,r}^{T}\}_{\bm n}$ is $d$-level $qr$-block GLT with symbol $[\kappa_{ij}]_{i,j=1}^q$ \cite[Theorem~4.24]{GLT20}.
		Here $e_i,e_j,e_\alpha$ are standard basis vectors in $\C^q$, $\C^{N(\bm n)}$, and $\C^r$, respectively.
		The permutation changes the order of the degrees of freedom from component first to grid point first.
		
		\item\label{glt:closure} \textbf{Closure under a.c.s. approximation.}
		A sequence $\{A_{\bm n}\}_{\bm n}$ is GLT with symbol $\kappa$ if and only if there exist GLT sequences $\{B_{\bm n,r}\}_{\bm n}\sim_{\mathrm{GLT}}\kappa_r$ such that $\{B_{\bm n,r}\}_{\bm n}\xrightarrow{\mathrm{a.c.s.}} \{A_{\bm n}\}_{\bm n}$ and $\kappa_r\to\kappa$ in measure \cite[Theorems~4.6 and~4.22]{GLT20}.
		
		\item\label{glt:isometry} \textbf{Matrix--symbol isometry.}
		For GLT sequences with symbols $\kappa$ and $\xi$,
		\[
		d_\psi(\{A_{\bm n}\},\{B_{\bm n}\})
		=d_\psi^{\mathrm{meas}}(\kappa,\xi)
		=\frac1{s|\Omega|}\int_\Omega
		\sum_{j=1}^s\psi(\sigma_j(\kappa(z)-\xi(z)))\,\mathrm dz.
		\]
		Thus a.c.s. convergence among GLT sequences is equivalent to convergence 	in measure of their symbols \cite[Remark~4.15 and Theorem~4.19]{GLT20}.
		
		\item\label{glt:approximation} \textbf{Elementary approximants.}
		Every GLT sequence with symbol $\kappa$ admits approximants
		\[
		B_{\bm n,r}=\sum_{j=1}^{J_r}D_{\bm n}(a_{j,r}I_s)T_{\bm n}(f_{j,r}),
		\qquad \kappa_r(x,\bm\vartheta)
		=\sum_{j=1}^{J_r}a_{j,r}(x)f_{j,r}(\bm\vartheta),
		\]
		where $J_r<\infty$, $a_{j,r}\in C^\infty([0,1]^d)$, and each $f_{j,r}$ is a monomial $e^{\mathrm i\ell\cdot\bm\vartheta}E_{\alpha\beta}$ with $\ell\in\mathbb Z^d$ and $E_{\alpha\beta}$ an $s\times s$ matrix unit. The approximants satisfy $\{B_{\bm n,r}\}_{\bm n}\xrightarrow{\mathrm{a.c.s.}} \{A_{\bm n}\}_{\bm n}$ and $\kappa_r\to\kappa$ almost everywhere \cite[Theorem~4.20]{GLT20}.
	\end{enumerate}
	The Karcher closure argument will use \ref{acs:schatten} to approximate the matrix means, \ref{glt:algebra}--\ref{glt:functional-calculus} to identify the symbols of fixed Richardson iterates, and \ref{glt:closure} to pass to their limit. The later squeezing argument uses \ref{glt:isometry} to compare ordered approximants even when their norms are not uniformly bounded.

	\subsection{The Karcher mean}
	\label{subsec:karcher-preliminaries}
	Let $\mathbb P_q$ be the cone of positive definite Hermitian matrices of order $q$. Its affine-invariant Riemannian distance is
	\begin{equation}
		\delta(A,B)=\bigl\|\log(A^{-1/2}BA^{-1/2})\bigr\|_2,
		\qquad A,B\in\mathbb P_q.
		\label{eq:riemannian-distance}
	\end{equation}
	The positive definite cone, with this metric, is a complete simply connected Riemannian manifold of nonpositive sectional curvature. Consequently, the sum of squared distances has a unique minimizer. The equal-weight \emph{Karcher mean}, also called the Riemannian least-squares mean, is
	\begin{equation}
		\label{eq:karcher-variational}
		\Karcher(B_1,\ldots,B_k)
		= \underset{X\in\mathbb P_q}{\operatorname{argmin}}\, \frac1{2k}\sum_{i=1}^k\delta(X,B_i)^2.
	\end{equation}
	This construction extends the midpoint interpretation of the two-matrix
	geometric mean; see \cite{BH06,LL11,Moakher05}. In particular,
	\[
	\Karcher(A,B) = A\#B =A^{1/2}(A^{-1/2}BA^{-1/2})^{1/2}A^{1/2}.
	\]
	The first-order condition for \eqref{eq:karcher-variational} is the logarithmic matrix equation in \ref{km:characterization}. For the properties below we use the local labels K1--K11, which follow from the variational definition. For their relation to the classical axioms for matrix geometric means, see \cite{ALM04,LL11,LL14}. Those axioms alone do not uniquely specify the Karcher mean.
	\begin{proposition}
		\label{prop:karcher-properties}
		Let $\mathbb P_q$ be the cone of positive definite Hermitian matrices of size $q$.
		Given $B_1,\ldots,B_k,C_1,\ldots,C_k\in\mathbb P_q$, the equal-weight mean $\Karcher(B_1,\ldots,B_k)$ satisfies:
		\begin{enumerate}[label=\textup{(K\arabic*)},ref=K\arabic*,leftmargin=*,itemsep=0.4em]
			\item \label{km:characterization}
			\textbf{Existence and uniqueness.}
			The mean is the unique matrix $Y\in\mathbb P_q$ satisfying
			\begin{equation}
				\sum_{i=1}^k\log(Y^{-1/2}B_iY^{-1/2})=O.
				\label{eq:glt-karcher-characterization}
			\end{equation}
			This is the first-order condition for \eqref{eq:karcher-variational}; see also 
			\cite[Theorem~2.9]{LP18}.
			
			\item \label{km:monotonicity}
			\textbf{Monotonicity in the L\"owner order.}
			If $B_i\preceq C_i$ for every $i$, then
			\[
			\Karcher(B_1,\ldots,B_k)\preceq\Karcher(C_1,\ldots,C_k).
			\]
			
			\item \label{km:congruence}
			\textbf{Congruence invariance.}
			For every invertible matrix $M\in\C^{q\times q}$,
			\begin{equation}
				\Karcher(M^*B_1M,\ldots,M^*B_kM)
				=M^*\Karcher(B_1,\ldots,B_k)M.
				\label{eq:glt-congruence}
			\end{equation}
			The matrix $M$ need not be unitary.
			
			\item \label{km:self-duality}
			\textbf{Inverse property.}
			\begin{equation}
				\Karcher(B_1^{-1},\ldots,B_k^{-1})
				=\Karcher(B_1,\ldots,B_k)^{-1}.
				\label{eq:glt-self-duality}
			\end{equation}
			
			\item \label{km:continuity}
			\textbf{Joint continuity.}
			The Karcher mean map $\Karcher:\mathbb P_q^k\to\mathbb P_q$ is continuous in the spectral norm topology. Thus, if $B_{i,r}\to B_i\succ O$ for every $i$, then $\Karcher(B_{1,r},\ldots,B_{k,r})\to\Karcher(B_1,\ldots,B_k)$.
			
			\item \label{km:commuting}
			\textbf{The commuting case.}
			If the $B_i$ commute pairwise, the Karcher mean has the following geometric-mean formula:
			\begin{equation}
				\Karcher(B_1,\ldots,B_k)
				=\exp\!\left(\frac1k\sum_{i=1}^k\log B_i\right)
				=(B_1\cdots B_k)^{1/k}.
				\label{eq:glt-commuting-mean}
			\end{equation}
			In particular, $\Karcher(cI,\ldots,cI)=cI$ for $c>0$.
			\item \label{km:symmetry}
			\textbf{Permutation invariance.} For every permutation $\pi$,
			\[
			\Karcher(B_{\pi(1)},\ldots,B_{\pi(k)})=\Karcher(B_1,\ldots,B_k).
			\]
			
			\item \label{km:homogeneity}
			\textbf{Joint homogeneity.} For $c_1,\ldots,c_k>0$,
			\[
			\Karcher(c_1B_1,\ldots,c_kB_k)
			=(c_1\cdots c_k)^{1/k}\Karcher(B_1,\ldots,B_k).
			\]
			
			\item \label{km:concavity}
			\textbf{Joint concavity.} If $0\le t\le1$, then
			\[
			\Karcher(tB_1+(1-t)C_1,\ldots,tB_k+(1-t)C_k)
			\succeq t\Karcher(B_1,\ldots,B_k)+(1-t)\Karcher(C_1,\ldots,C_k).
			\]
			
			\item \label{km:determinant}
			\textbf{Determinant identity.}
			\[
			\det\Karcher(B_1,\ldots,B_k)
			=\left(\prod_{i=1}^k\det B_i\right)^{1/k}.
			\]
			
			\item \label{km:agh}
			\textbf{Harmonic and arithmetic bounds.}
			\[
			\left(\frac1k\sum_{i=1}^k B_i^{-1}\right)^{-1}
			\preceq\Karcher(B_1,\ldots,B_k)
			\preceq\frac1k\sum_{i=1}^k B_i.
			\]
		\end{enumerate}
	\end{proposition}
	
	Properties K1--K6 follow from \cite[Theorems~6.5 and~6.8, p.~19]{LL14}.
	The former gives \ref{km:characterization}; the latter gives \ref{km:monotonicity}, \ref{km:congruence}, \ref{km:self-duality}, \ref{km:continuity}, and \ref{km:commuting} as properties (P4), (P6), (P8), (P5), and (P1), respectively. For \ref{km:continuity}, we also use the fact that the metric
	\[
	d_T(B,C)=\|\log(B^{-1/2}CB^{-1/2})\|
	\]
	induces the norm topology on $\mathbb P_q$ \cite[Section 2, p.~4]{LL14}.
	We also recall three basic matrix properties: congruence preserves the L\"owner order; inversion reverses it, $O\prec B\preceq C\Rightarrow C^{-1}\preceq B^{-1}$; and
	\[
	(M^*BM)^{-1}=M^{-1}B^{-1}(M^*)^{-1}.
	\]
	By \ref{km:monotonicity} and \ref{km:commuting}, $aI\preceq B_i\preceq bI$ implies $aI\preceq\Karcher(B_1,\ldots,B_k)\preceq bI$.
	
	Permutation invariance, joint homogeneity, the determinant identity, joint concavity, and the harmonic--arithmetic bounds are reported and discussed in \cite{LL11,ALM04}.
	
	The closure proof uses monotonicity to compare truncated inputs and self-duality to transfer these comparisons through inversion. Continuity is also used in our arguments at fixed block size.

	\section{Dimension-independent Schatten estimates}\label{sec:schatten-main}
	The following theorem is both the main result of this work and the bottleneck needed to achieve the GLT closure of the Karcher mean. This section is devoted to presenting this result, whose proof is preceded by some preparatory lemmas.
	\begin{theorem}\label{th:task}
		Let $A_{\bm n}^{(i)},X_{\bm n,0} \in \C^{sN(\bm n)\times sN(\bm n)}$, for $i=1,\ldots,k$, be HPD matrices. Suppose that there are constants $0<a\le b<\infty$, independent of $\bm n$, such that
		\begin{equation} \label{eq:bounds} 
			aI\preceq A_{\bm n}^{(i)}\preceq bI, \qquad	aI\preceq X_{\bm n,0}\preceq bI.
		\end{equation}
		Set
		\begin{displaymath}
			R=\log(b/a),
			\qquad
			0<\hat{\theta}\le \frac{1}{2k(1+2R)},
		\end{displaymath}
		and define
		\begin{align}
			X_{\bm n,\nu+1}
			={}&X_{\bm n,\nu}
			-\hat{\theta} X_{\bm n,\nu}^{1/2}
			\left[
			\sum_{i=1}^{k}
			\log\!\left(
			X_{\bm n,\nu}^{1/2}
			(A_{\bm n}^{(i)})^{-1}
			X_{\bm n,\nu}^{1/2}
			\right)
			\right]
			X_{\bm n,\nu}^{1/2}.
			\tag{R}\label{eq:R}
		\end{align}
		Let $G_{\bm n}=\Karcher\!\left(A_{\bm n}^{(1)},\ldots,A_{\bm n}^{(k)}\right)$. Then, for every $\bm n$ and every $\nu\ge0$, we have the spectral bounds
		\begin{equation}
			aI\preceq X_{\bm n,\nu}\preceq bI, \qquad aI\preceq G_{\bm n}\preceq bI.
			\label{eq:interval}
		\end{equation}
		Moreover, the inequality
		\begin{equation}
			\frac{\|X_{\bm n,\nu}-G_{\bm n}\|_2^2}{sN(\bm n)} \leq b^2R^2\rho^\nu
			\label{eq:frob}
		\end{equation}
		holds {with $\rho=1-k\hat{\theta}$}, and, for every $1\leq p<\infty$,
		\begin{equation}
			\frac{\norm{X_{\bm n,\nu}-G_{\bm n}}_p}{(sN(\bm n))^{1/p}} \leq bR\,\rho^{\nu/\max\{2,p\}}.
			\label{eq:schatten}
		\end{equation}
	\end{theorem}
	
	\noindent {We recall that the formula \eqref{eq:R} is in fact the iteration in \cite[equation~(4.20), p.~51]{TH}. Moreover, given a set $ \{ A_{\bm n}^{(i)}\}_{i=1}^{k} $ of HPD matrices, the initialization point $ X_{\bm n}^{(0)} $ for the iteration \eqref{eq:R} can be set as the arithmetic mean
		\begin{displaymath}
			X_{\bm n,0}=\frac1k\sum_{i=1}^{k}A_{\bm n}^{(i)}
		\end{displaymath}
		respecting the spectral bounds \eqref{eq:bounds}.}
	
	The result in \cite[Theorem~2, p.~1703]{BI13} is local and does not {imply by itself Theorem \ref{th:task}. In fact, in the the proof we use} the auxiliary estimates below to control the actual step in \eqref{eq:R}. By exploiting interval invariance, we keep all later estimates on the same spectral interval. We then obtain the energy decay and the metric comparison, which we then convert to the stated Frobenius and Schatten bounds.
	
	{Unless specified, we will write $A_i$ to indicate $A_{\bm n}^{(i)}$ for the sake of notational simplicity.}

	\subsection{Auxiliary estimates}
	{We provide here} the estimates used in the proof of Theorem \ref{th:task}; notably, spectral interval invariance, control of the objective gap and decrease along the affine step. The metric comparison of Lemma \ref{lem:distance-comparison} then gives the Frobenius estimate. Throughout, we retain the input and step-size hypotheses of Theorem~\ref{th:task}.
	
	
	The first task is to keep every iterate in the prescribed spectral interval. This ensures that the next step is defined and that the same spectral bounds remain available throughout the argument. For $X\succ O$, define the Hermitian matrix
	\[
	V(X)=\frac1k\sum_{i=1}^{k}\log\!\left(X^{-1/2}A_iX^{-1/2}\right).
	\]
	Since {$X^{1/2}A_i^{-1}X^{1/2}=\left(X^{-1/2}A_i X^{-1/2}\right)^{-1}$,} the identity $\log(C^{-1})=-\log C$ gives the step:
	\begin{equation}
		X_+
		=X+\alpha X^{1/2}V(X)X^{1/2}
		=X^{1/2}\bigl(I+\alpha V(X)\bigr)X^{1/2},
		\label{eq:update}
	\end{equation}
	{with $\alpha=k\hat{\theta}$.}
	
	\begin{lemma}\label{lem:interval-invariance}
		Let $aI\preceq X\preceq bI$. Under the hypotheses of Theorem~\ref{th:task}, \eqref{eq:update} satisfies
		\begin{equation}
			aI\preceq X_+\preceq bI;
			\label{eq:invariant}
		\end{equation}
		Moreover, the update direction satisfies
		\begin{equation}
			-RI\preceq V(X)\preceq RI,
			\qquad
			\|V(X)\|\le R.
			\label{eq:Vbound}
		\end{equation}
		The geometric mean also satisfies $aI\preceq G\preceq bI$ and the bounds in \eqref{eq:interval} hold for every iterate $X_{ \nu}$.
	\end{lemma}
	
	\begin{proof}
		{Since $aI\preceq X\preceq bI$ by hypothesis, then $aX^{-1}\preceq X^{-1/2}A_iX^{-1/2}\preceq bX^{-1}$.} The operator monotonicity of the logarithm gives
		\begin{equation}
			(\log a)I-\log X
			\preceq V(X)
			\preceq (\log b)I-\log X;
			\label{eq:logorder}
		\end{equation}
		indeed, it follows directly from the resolvent representation: if
		$O\prec A\preceq B$, then
		\[
		\log B-\log A
		=\int_0^\infty\Bigl((A+tI)^{-1}-(B+tI)^{-1}\Bigr)\,\dd t
		\succeq O.
		\]
		For $c>0$, set $h_c(x)=x+\alpha x\log(c/x)$. Applying congruence by $X^{1/2}$ to
		\eqref{eq:logorder} gives
		\begin{equation}
			h_a(X)\preceq X_+\preceq h_b(X).
			\label{eq:horder}
		\end{equation}
		On the interval $[a,b]$, since $\alpha\le 1/[2(1+2R)]\le1/(1+R)$,
		\[
		h_a'(x)=1+\alpha\bigl(\log(a/x)-1\bigr)
		\ge 1-\alpha(1+R)\ge0;
		\]
		hence $h_a(x)\ge h_a(a)=a$.
		Similarly,
		\[
		h_b'(x)=1+\alpha\bigl(\log(b/x)-1\bigr)
		\ge1-\alpha\ge0,
		\]
		so $h_b(x)\le h_b(b)=b$. Computations for the single matrix $X$ in \eqref{eq:horder} therefore gives \eqref{eq:invariant}; by induction, all iterates are positive definite and remain in the prescribed interval. Also,
		\[
		e^{-R}I\preceq X^{-1/2}A_iX^{-1/2}\preceq e^RI,
		\]
		which leads to \eqref{eq:Vbound}.
		
		\noindent By Definition \ref{eq:karcher-variational}, the Karcher mean is the unique positive definite solution of
		\begin{equation}
			V(G)=O.
			\label{eq:karcher-eq}
		\end{equation}
		Existence and uniqueness follow from \cite[Theorem~2.9]{LP18}, applied to the finitely supported measure $k^{-1}\sum_i\delta_{A_i}$. The considered measure satisfies {the assumptions} of the cited result: its Hilbert-space operator setting includes $\C^m$, and hence the complex Hermitian matrices considered here.
		
		Applying the logarithmic bounds to $G$, without assuming in advance that $G$ lies in the interval, gives
		\[
		(\log a)I-\log G\preceq O\preceq(\log b)I-\log G;
		\]
		This chain of inequalities concerns functions of a single Hermitian matrix, that $aI\preceq G\preceq bI$. This proves \eqref{eq:interval}.
	\end{proof}
	
	\noindent {As a subsequent step, we} recall the distance
	\[
	\delta(X,A)=\left\|\log\!\left(X^{-1/2}AX^{-1/2}\right)\right\|_2
	\]
	and the normalized objective {function}
	\begin{equation}
		f(X)=\frac1{2k}\sum_{i=1}^{k}\delta(X,A_i)^2.
		\label{eq:fdef}
	\end{equation}
	
	
	\begin{definition}[Spectral derivative]
		Let $F:J \to \mathbb{R}$ be smooth on an open interval $J\subset\mathbb R$, and let $\Phi(X)=F(X)$ be the spectral function for {a Hermitian matrix} $X$ with spectrum in $J$. The \emph{spectral derivative}, also called the Fr\'echet derivative of the spectral function, is the linear map
		\begin{displaymath}
			D\Phi(X)\left[H\right] = Q\left(F^{\left[1\right]}(\Lambda) \circ \left(Q^* H Q\right) \right)Q^*,
		\end{displaymath}
		where $\circ$ denotes the entrywise Hadamard product, $\Lambda$ is diagonal and satisfies $X = Q\Lambda Q^*$ with $Q$ unitary, and $[F^{[1]}(\Lambda)]_{j\ell}=F^{[1]}(\lambda_j,\lambda_\ell)$ for the diagonal entries $\lambda_j$ of $\Lambda$ {defined as}
		\begin{equation}\label{frechet}
			F^{\left[1\right]}(x,y) = 
			\begin{cases}
				\frac{F(x)-F(y)}{x-y} & x\neq y,\\
				F'(x) & x=y.
			\end{cases}
		\end{equation}
	\end{definition}
	
	\begin{lemma}
		Let
		\[
		f_A(X)=\frac12\delta(X,A)^2,
		\qquad
		U=X^{-1/2}HX^{-1/2},
		\qquad
		B=\log\!\left(X^{-1/2}AX^{-1/2}\right),
		\]
		with $H$ Hermitian. Choose an orthonormal eigenbasis of $B$, write $\lambda_1,\ldots,\lambda_m$ for its eigenvalues, and let $u_{j\ell}$ be the entries of $U$ in this basis. Set
		\[
		\chi(t)=
		\begin{cases}
			t\coth t, & t\ne0,\\
			1, & t=0.
		\end{cases}
		\]
		{Then}
		\begin{equation}\label{eq:first-derivative}
			Df_A(X)[H]=-\tr(BU),
		\end{equation}
		and the ordinary affine second derivative is
		\begin{equation}\label{eq:affine-hessian}
			D^2f_A(X)[H,H]
			=\sum_{j,\ell=1}^{m}
			\left[
			\chi\!\left(\frac{\lambda_j-\lambda_\ell}{2}\right)
			+\frac{\lambda_j+\lambda_\ell}{2}
			\right]|u_{j\ell}|^2.
		\end{equation}
		Along the curve $\gamma(t)=X^{1/2}e^{tU}X^{1/2}$, we have
		\begin{equation}\label{eq:geo-hessian}
			\left.\frac{\dd^2}{\dd t^2}f_A(\gamma(t))\right|_{t=0}
			=\sum_{j,\ell=1}^{m}
			\chi\!\left(\frac{\lambda_j-\lambda_\ell}{2}\right)|u_{j\ell}|^2.
		\end{equation}
	\end{lemma}
	
	\begin{proof}
		The quantity $\delta(X,A)^2$ is the sum of the squared logarithms of the eigenvalues of $X^{-1}A$. As a consequence, 
		$\delta(X,A)^2$ is therefore invariant {under invertible congruence and symmetric in both $X$ and $A$}, since interchanging the two matrices inverts these eigenvalues. In particular, $f_A(X+tH)=f_{e^B}(I+tU)$. We may thus reduce to
		$X=I$, $A=e^B$, and then diagonalize $B$ by a unitary matrix.
		
		Set $q(x)=\tfrac12(\log x)^2$. Along the affine curve,
		\begin{equation}\label{eq:faitu}
			f_{e^B}(I+tU)=\tr q\!\left(e^{-B/2}(I+tU)e^{-B/2}\right).
		\end{equation}
		Now let $Y(t) = e^{-B/2}(I+tU)e^{-B/2}$. Then $Y(0)=e^{-B}$ and $Y'(t) = e^{-B/2}Ue^{-B/2} = Y'(0)$. Hence
		\begin{equation}
			\begin{aligned}
				\left.\frac{\dd}{\dd t}f_{e^B}(I+tU)\right|_{t=0} &= \left.\frac{\dd}{\dd t}\tr q(Y(t))\right|_{t=0} = \tr  \left(q'\!\left(Y(0)\right)Y'(0)\right)\\
				&= \tr \left(\log(e^{-B}) \left(e^{-B}\right)^{-1}  e^{-B/2}Ue^{-B/2}\right)\\
				&= \tr \left(-B e^{B} e^{-B/2}Ue^{-B/2}\right)\\
				&=  - \tr \left(e^{-B/2} B e^{B/2}U\right) \mathrel{=} - \tr \left(BU\right),
			\end{aligned}
		\end{equation}
		The first spectral derivative therefore gives \eqref{eq:first-derivative}. For the second derivative, we proceed as follows. The definition of $Y(t)$ gives $Y''=0$. {Setting $d_j=e^{-\lambda_j}$ and $K=Y'(0)$, by} \eqref{frechet},
		\begin{equation}
			D^2f_A(X)[H,H] \mathrel{=} \sum_{i,j} q'^{\left[1\right]}(d_i,d_j)\lvert K_{ij}\rvert^2,
		\end{equation}
		where $K_{ij}$ is {the entry $(K)_{ij}$ with respect to the eigenbasis} of $B=\text{diag}(\lambda_1,\dots,\lambda_m)$, namely
		\begin{equation}
			K_{ij} = e^{-\lambda_i/2}u_{ij} e^{-\lambda_j/2},
		\end{equation}
		and
		\begin{equation}
			q'^{\left[1\right]}(x,y) = 
			\begin{cases}
				\frac{q'(x)-q'(y)}{x-y} &x\neq y,\\
				q''(x) & x=y.
			\end{cases}
		\end{equation}
		Thus
		\begin{equation}
			e^{-B} = \text{diag}\left(e^{-\lambda_1},\dots, e^{-\lambda_m}\right) = Y(0).
		\end{equation}
		If $d_j = e^{-\lambda_j}$ and $q'(x)=x^{-1}\log x$, then $q'(d_j) = -\lambda_je^{\lambda_j}$, and for $\lambda_i\ne\lambda_j$
		\begin{equation}
			q'^{\left[1\right]}(d_i,d_j) = \frac{-\lambda_i e^{\lambda_i}  + \lambda_j e^{\lambda_j}}{e^{-\lambda_i} - e^{-\lambda_j}}.
		\end{equation}
		Multiplication by $\lvert K_{ij}\rvert^2$ leads to the equality
		\begin{equation}
			\lvert K_{ij}\rvert^2 q'^{\left[1\right]}(d_i,d_j) = e^{-\lambda_i}e^{-\lambda_j} \lvert u_{ij}\rvert^2 \frac{-\lambda_i e^{\lambda_i}  + \lambda_j e^{\lambda_j}}{e^{-\lambda_i} - e^{-\lambda_j}}.
		\end{equation}
		Since {$e^{-\lambda_i} - e^{-\lambda_j} = e^{-\lambda_i - \lambda_j}\left(e^{\lambda_j} - e^{\lambda_i}\right)$,} we obtain
		\begin{equation}\label{eq:der_sec_frac}
			\lvert u_{ij}\rvert^2  e^{-\lambda_i - \lambda_j} \frac{-\lambda_i e^{\lambda_i}  + \lambda_j e^{\lambda_j}}{e^{-\lambda_i - \lambda_j}\left(e^{\lambda_j} - e^{\lambda_i}\right)} = \lvert u_{ij}\rvert^2  \frac{\lambda_i e^{\lambda_i}-\lambda_j e^{\lambda_j}}
			{e^{\lambda_i}-e^{\lambda_j}}.
		\end{equation}
		Setting
		\begin{equation}
			\mu = \frac{\lambda_i+\lambda_j}{2} \qquad \tau = \frac{\lambda_i-\lambda_j}{2} \qquad \Longrightarrow \qquad\lambda_i=\mu+\tau \qquad \lambda_j = \mu-\tau,
		\end{equation}
		we can rewrite the numerator and denominator in \eqref{eq:der_sec_frac} as
		\begin{equation}
			\begin{aligned}
				\lambda_i e^{\lambda_i}-\lambda_j e^{\lambda_j} &= (\mu+\tau)e^{\mu+\tau} - (\mu-\tau)e^{\mu-\tau} = 2e^\mu\left(\mu \sinh(\tau) + \tau\cosh(\tau)\right),\\
				e^{\lambda_i}-e^{\lambda_j} &= e^{\mu+\tau} - e^{\mu-\tau} = 2e^\mu \sinh(\tau).
			\end{aligned}
		\end{equation}
		Thus, for distinct eigenvalues, we obtain
		\begin{equation}
			\frac{2e^\mu\left(\mu \sinh(\tau) + \tau\cosh(\tau)\right)}{2e^\mu \sinh(\tau)} = \mu + \tau\coth(\tau) = \frac{\lambda_i+\lambda_j}{2}
			+\chi\!\left(\frac{\lambda_i-\lambda_j}{2}\right).
		\end{equation}
		For coincident eigenvalues, on the other hand, the continuous extension is $1+\lambda_j$, since $\chi(0)=1$. This gives \eqref{eq:affine-hessian}, including the case of repeated eigenvalues.
		
		For the curve $e^{tU}$, the {second derivative} at zero is $U^2$. The additional term from the chain rule is
		\[
		Df_{e^B}(I)[U^2]
		=-\tr(BU^2)
		=-\frac12\sum_{j,\ell}(\lambda_j+\lambda_\ell)|u_{j\ell}|^2.
		\]
		This cancels the second term in \eqref{eq:affine-hessian}, proving \eqref{eq:geo-hessian}.
	\end{proof}
	
	\noindent The elementary inequality
	\begin{equation}
		1\le\chi(t)\le1+|t|
		\label{eq:psi-bound}
	\end{equation}
	gives the two estimates needed below. Indeed, for $t\ge0$,
	\[
	t\coth t=t+\frac{2t}{e^{2t}-1}\le t+1;
	\]
	the lower bound follows from $\tanh t\le t$ and the evenness of $\chi$. Equation \eqref{eq:geo-hessian} gives the global lower bound along geodesics
	\begin{equation}
		\left.\frac{\dd^2}{\dd t^2}f_A(\gamma(t))\right|_{t=0}
		\ge\|U\|_2^2.
		\label{eq:geo-lower}
	\end{equation}
	If, in addition, $aI\preceq X,A\preceq bI$, then $|\lambda_j|\le R$. Equations \eqref{eq:affine-hessian} and \eqref{eq:psi-bound} give
	\[
	D^2f_A(X)[H,H]
	\le(1+2R)\|X^{-1/2}HX^{-1/2}\|_2^2.
	\]
	Averaging over the inputs yields
	\begin{equation}
		D^2f(X)[H,H]
		\le (1+2R)\|X^{-1/2}HX^{-1/2}\|_2^2,
		\qquad aI\preceq X\preceq bI.
		\label{eq:affine-upper}
	\end{equation}
	The term $(\lambda_j+\lambda_\ell)/2$ in \eqref{eq:affine-hessian} is essential for
	estimating the additive step. Omitting it would give the estimate for the exponential curve.
	
	We now compare the objective gap with the distance to the mean and with the
	size of the update direction. These comparisons connect the decrease of
	the objective to an error estimate without assuming that the current
	matrix is close to the mean.
	
	\begin{lemma}\label{lem:objective-gap}
		For every {$X \succ O$}, the inequality below hold uniformly
		\begin{equation}\label{eq:two-global}
			\frac12\delta(X,G)^2 \le f(X)-f(G) \le\frac12\norm{V(X)}_2^2.
		\end{equation}
	\end{lemma}
	
	\begin{proof}
		For arbitrary $X,Y\succ O$, {define}
		\begin{equation}
			Z=\log(X^{-1/2}YX^{-1/2}), \qquad \gamma(t)=X^{1/2}e^{tZ}X^{1/2}.
		\end{equation}
		Then $\gamma(0) = X$ and $\gamma(1) = Y$. With $M_t=X^{1/2}e^{tZ/2}$, we have
		\begin{enumerate}
			\item $\gamma(t)=M_tM_t^*$;
			
			\item $\gamma'(t) = X^{1/2} Z e^{tZ} X^{1/2}= M_t Z M_t^*$, since $Z$ commutes with $e^{tZ}$;
			
			\item $\gamma^{-1/2}M_t$ is unitary, since $(\gamma^{-1/2}M_t)(\gamma^{-1/2}M_t)^*=\gamma^{-1/2}(M_tM_t^*)\gamma^{-1/2}=I$;
			
			\item $\gamma''(t) = \gamma'(t)\gamma^{-1}(t)\gamma'(t)$.
		\end{enumerate}
		Thus \eqref{eq:geo-lower} applies at every point of the curve, with $\|\gamma(t)^{-1/2}\gamma'(t)\gamma(t)^{-1/2}\|_2=\|Z\|_2$ because the normalized velocity is unitarily similar to $Z$. Integrating its averaged form and using \eqref{eq:first-derivative} gives
		\begin{equation}\label{eq:strongconvex}
			f(Y)\ge f(X)-\tr(V(X)Z)+\frac12\norm{Z}_2^2.
		\end{equation}
		Since
		\begin{equation}
			\norm{Z}_2^2 = \norm{\log\left( X^{-1/2}Y X^{-1/2} \right)}_2^2 = \delta(X,Y)^2 = \delta(Y,X)^2,
		\end{equation}
		taking $X=G$ and using $V(G)=O$ gives
		\begin{equation}\label{eq:error-lower}
			f(Y)-f(G)\ge\frac12\delta(Y,G)^2.
		\end{equation}
		In particular, the solution of the Karcher equation is the global minimizer of $f$.
		
		Next, set $Y=G$ in \eqref{eq:strongconvex}. The Cauchy--Schwarz inequality for the Frobenius inner product gives
		\begin{equation}\label{eq:error-upper}
			\begin{aligned}
				f(X)-f(G) & \le \tr(V(X)Z)-\frac12\norm{Z}_2^2\\
				&\le \norm{V(X)}_2\norm{Z}_2-\frac12\norm{Z}_2^2\\
				&\le \frac12\norm{V(X)}_2^2.
			\end{aligned}
		\end{equation}
		Combining the two bounds gives \eqref{eq:two-global}.
	\end{proof}
	
	The affine second-derivative bound controls the remainder in Taylor's formula along the actual update segment; the step-size restriction makes this remainder small enough to retain a decrease in the {objective function \eqref{eq:fdef}}. Combining that decrease with the preceding objective-gap bound gives a geometric decay estimate, as detailed in the next lemma.
	
	\begin{lemma}\label{lem:affine-descent}
		Let $ X $ be HPD such that $X\in[aI,bI]$ with respect to the L\"owner order $ \succ $; then additive update satisfies
		\begin{equation}\label{eq:descent}
			f(X_+)\le f(X)-\frac\alpha2\norm{V(X)}_2^2.
		\end{equation}
		Consequently, the iterates satisfy
		\begin{equation}\label{eq:energy}
			f(X_\nu)-f(G) \le\rho^\nu\bigl(f(X_0)-f(G)\bigr).
		\end{equation}
	\end{lemma}
	
	\begin{proof}
		Let {$aI\preceq X\preceq bI$}. Write $V=V(X)$ and set
		\begin{displaymath}
			H=\alpha X^{1/2}VX^{1/2}, \qquad X_t=X+tH, \qquad 0\le t\le1.
		\end{displaymath}
		The invariance proved above places both $X$ and $X+H$ in $[aI,bI]$. The entire affine segment $X_t$ therefore remains in this interval. By \eqref{eq:first-derivative},
		\begin{equation}\label{eq:descent-first}
			Df(X)[H]=-\alpha\norm{V}_2^2.
		\end{equation}
		Also, $X_t=X^{1/2}(I+t\alpha V)X^{1/2}$. Using cyclicity of the trace and the fact that $V$ commutes with $C_t:=I+t\alpha V$,
		\begin{displaymath}
			\norm{X_t^{-1/2}HX_t^{-1/2}}_2^2 = \tr\!\left( \left(X_t^{-1} H\right)^2\right).
		\end{displaymath}
		Now
		\begin{displaymath}
			X_t^{-1}H = X^{-1/2} \left(I + t\alpha V\right)^{-1} X^{-1/2} \alpha X^{1/2} V X^{1/2} = \alpha X^{-1/2} C_t^{-1} V X^{1/2};
		\end{displaymath}
		hence
		\begin{displaymath}
			\left(X_t^{-1}H\right)^2 = \alpha^2 X^{-1/2} C_t^{-1} V C_t^{-1} V X^{1/2},
		\end{displaymath}
		and, again using cyclicity of the trace and commutativity of $V$ with $C_t$,
		\begin{equation}\label{eq:normalized-step}
			\begin{aligned}
				\tr\!\left(\left(X_t^{-1}H\right)^2\right) &= \alpha^2 \tr\!\left(\left(C_t^{-1}V\right)^2\right) = \alpha^2\tr\!\left(V^2 C_t^{-2}\right)\\
				&= \alpha^2 \tr\!\left(V^2 \left(I + t\alpha V\right)^{-2}\right)\\
				&\leq \frac{\alpha^2}{(1-\alpha R)^2}\norm{V}_2^2.
			\end{aligned}
		\end{equation}
		Indeed,
		\begin{displaymath}
			\alpha R\le\frac{R}{2(1+2R)}<\frac14,
		\end{displaymath}
		and all denominators are positive.
		Taylor's formula with integral remainder along this affine segment, together with \eqref{eq:affine-upper}, gives
		\begin{equation}\label{eq:taylor}
			f(X+H) \le f(X)-\alpha\norm{V}_2^2 +\frac{(1+2R)\alpha^2}{2(1-\alpha R)^2}\norm{V}_2^2.
		\end{equation}
		The step-size restriction implies
		\begin{displaymath}
			\frac{\alpha (1+2R)}{(1-\alpha R)^2} \le\frac{1/2}{(3/4)^2} =\frac89<1.
		\end{displaymath}
		This proves \eqref{eq:descent}.
		Combining this estimate with \eqref{eq:error-upper}, we deduce
		\begin{equation}
			\begin{aligned}
				f(X_+)-f(G)
				&\le f(X)-f(G)-\frac\alpha2\norm{V(X)}_2^2\\
				&\le(1-\alpha)\bigl(f(X)-f(G)\bigr).
			\end{aligned}
		\end{equation}
		Iteration gives \eqref{eq:energy}.
	\end{proof}
	
	\noindent Note that normalization in \eqref{eq:fdef} accounts for the factor $k$: the step size for the averaged objective is $\alpha=k\hat{\theta}$, and the convergence factor is exactly $\rho=1-k\hat{\theta}$.
	
	The following comparison converts a Riemannian-distance bound into a
	Frobenius bound with the upper spectral bound $b$ as its constant.
	
	\begin{lemma}\label{lem:distance-comparison}
		Let $X,G\in[aI,bI]$. Then
		\begin{equation}\label{eq:euclidean-vs-riemannian}
			\norm{X-G}_2 \leq b \delta(X,G).
		\end{equation}
	\end{lemma}
	
	\begin{proof}
		Set
		\begin{displaymath}
			C=X^{-1/2}GX^{-1/2}, \qquad \gamma(t)=X^{1/2}C^tX^{1/2}, \qquad 0\le t\le1.
		\end{displaymath}
		The scalar inequality $u^t\le(1-t)+tu$, applied to $C$, gives {$\gamma(t)\preceq(1-t)X+tG\preceq bI$.}
		With $M_t=X^{1/2}C^{t/2}$, we have
		\begin{displaymath}
			\gamma'(t)=M_t(\log C)M_t^*, \qquad \norm{M_t}^2=\norm{\gamma(t)}\le b.
		\end{displaymath}
		Thus
		\begin{displaymath}
			\norm{\gamma'(t)}_2 \leq b\norm{\log C}_2 = b \delta(X,G).
		\end{displaymath}
		Integrating from $0$ to $1$ gives \eqref{eq:euclidean-vs-riemannian}.
	\end{proof}

	\subsection{Proof of the Schatten estimates}
	The energy estimate controls the objective gap but still depends on its initial value. We bound that value using the spectral interval, then obtain a Riemannian-distance estimate and compare it with the Frobenius norm on the same interval.
	
	\begin{proof}[Proof of Theorem~\ref{th:task}]
		By Lemma~\ref{lem:interval-invariance}, both the iterates and the mean remain in the spectral interval. Lemma~\ref{lem:affine-descent} gives the decay of the objective gap, and Lemma~\ref{lem:objective-gap} bounds the squared Riemannian distance by that gap. We first bound the initial objective gap.
		Since $X_0,A_i\in[aI,bI]$, all eigenvalues of $X_0^{-1/2}A_iX_0^{-1/2}$ lie in $[e^{-R},e^R]$. Hence
		\begin{displaymath}
			\delta(X_0,A_i)^2\le mR^2,
		\end{displaymath}
		and therefore
		\begin{equation}\label{eq:init}
			0\le f(X_0)-f(G)\le f(X_0)\le\frac{mR^2}{2}.
		\end{equation}
		Equations \eqref{eq:error-lower}, \eqref{eq:energy}, and \eqref{eq:init} give
		\begin{equation}\label{eq:delta-rate}
			\delta(X_\nu,G)^2\le mR^2\rho^\nu.
		\end{equation}
		Combining Lemma~\ref{lem:distance-comparison} with \eqref{eq:delta-rate} yields
		\begin{displaymath}
			\norm{X_\nu-G}_2^2 \leq b^2mR^2\rho^\nu.
		\end{displaymath}
		{Setting} $m=sN(\bm n)$ proves \eqref{eq:frob}, with the initial error bound requiring only that $X_0\in [aI,bI]$.
		The arithmetic mean is an admissible initial matrix because this interval is convex.
		
		For indices at most two, the normalized Frobenius estimate gives the Schatten estimate by comparison of vector norms; for larger indices, we also use the spectral-norm bound supplied by interval invariance.
		Let $E_\nu=X_\nu-G$; for $1\le p\le2$, the norm inequality for the vector of singular values gives
		\begin{equation}\label{eq:p-le-2}
			m^{-1/p}\norm{E_\nu}_p \leq m^{-1/2}\norm{E_\nu}_2 \leq bR\rho^{\nu/2}.
		\end{equation}
		For $p>2$, invariance of the interval gives
		\begin{displaymath}
			-(b-a)I\preceq E_\nu\preceq(b-a)I,
		\end{displaymath}
		so that
		\begin{equation}\label{eq:op-error}
			\norm{E_\nu}\leq b-a=b(1-e^{-R})\leq bR.
		\end{equation}
		As a consequence we obtain
		\begin{equation}
			\begin{aligned}
				\norm{E_\nu}_p^p &=\sum_{j=1}^{m}\sigma_j(E_\nu)^p\\
				&\leq \norm{E_\nu}^{p-2} \norm{E_\nu}_2^2\\
				&\leq (bR)^{p-2}b^2mR^2\rho^\nu\\
				&= m(bR)^p\rho^\nu,
			\end{aligned}
		\end{equation}
		and hence
		\begin{equation}\label{eq:p-gt-2}
			m^{-1/p}\norm{E_\nu}_p\le bR\rho^{\nu/p}.
		\end{equation}
		Thus Equations \eqref{eq:p-le-2} and \eqref{eq:p-gt-2} prove \eqref{eq:schatten}, completing the proof of Theorem~\ref{th:task}.
	\end{proof}

	\subsection{Some remarks}
	\paragraph{Local convergence is not sufficient.}
	The result in \cite[Theorem~2, p.~1703]{BI13} assumes that the initial matrix lies in a sufficiently small neighborhood of the mean $G$, with radius given by
	\begin{displaymath}
		\norm{G^{-1/2}X_0G^{-1/2}-I}_2.
	\end{displaymath}
	This alone does not imply an arbitrary initial matrix in $[aI,bI]$ or the global, dimension-independent remainder estimate needed here. The additional argument we need is the affine bound \eqref{eq:affine-upper}, followed by the affine-step calculation in Lemma~\ref{lem:affine-descent}.
	
	The statement following \cite[equation~(4.20), p.~51]{TH} should therefore not be read as convergence from every positive definite initial matrix with an arbitrary fixed step size; indeed, for scalar inputs $A_i=1$, any fixed $\hat{\theta}>0$ gives
	\begin{displaymath}
		x_+=x\bigl(1-k\hat{\theta}\log x\bigr),
	\end{displaymath}
	and the choice $x_0=e^{2/(k\hat{\theta})}$ gives $x_1=-x_0<0$. Theorem~\ref{th:task} excludes this failure by choosing the step size from an interval containing both the inputs and the initial matrix.

	\paragraph{Normalized approximation for GLT closure.}
	By \eqref{eq:schatten}, for every $1\le p<\infty$,
	\begin{displaymath}
		\frac{\norm{X_{\bm n,\nu}-G_{\bm n}}_p^p}{sN(\bm n)} \leq (bR)^p\rho^{\nu p/\max\{2,p\}};
	\end{displaymath}
	since $0<\rho<1$, it follows that
	\begin{equation}
		\lim_{\nu\to\infty}\limsup_{\bm n\to\infty} \frac{\norm{X_{\bm n,\nu}-G_{\bm n}}_p^p}{sN(\bm n)} =0, \qquad 1\leq p<\infty.
		\label{eq:acs}
	\end{equation}
	The bound is uniform in the matrix size, so the iteration count can be chosen independently of that dimension. Equation \eqref{eq:acs} gives us the Schatten estimate for the a.c.s. argument in Section~\ref{sec:glt-karcher}. To identify the GLT symbol,
	we must also show that each fixed iterate is GLT and that its symbol converges.

	\paragraph{Quantitative estimate.}
	The preceding argument gives the stated estimates for the additive Richardson
	iteration on complex Hermitian matrices. The step size
	depends only on $a,b,k$, while the normalized constants and convergence factor are
	independent of $\bm n$ and $m=sN(\bm n)$. Uniform spectral bounds are used to prove
	invariance of the interval, control the affine second derivative and bound the initial error.

	\section{GLT closure and the symbol of the Karcher mean}\label{sec:glt-karcher}
	The Schatten estimates allow us to identify the GLT symbol of the Karcher mean for a fixed finite collection of input sequences. We first treat uniform spectral bounds, then remove them by squeezing the mean between two GLT approximations.
	We use the notation and Karcher properties from Section~\ref{sec:preliminaries}.

	\subsection{GLT closure under uniform spectral bounds}
	The proof combines the Schatten criterion \ref{acs:schatten} with the GLT algebra and functional-calculus rules \ref{glt:algebra}--\ref{glt:functional-calculus}.
	We first assume uniform positive lower and finite upper spectral bounds on all input sequences.
	
	\begin{theorem}\label{thm:glt-bounded}
		Let $A_{\bm n}^{(i)}$ be HPD matrices such that, for $i=1,\ldots,k$,
		\[
		\{A_{\bm n}^{(i)}\}_{\bm n}\sim_{\mathrm{GLT}}\kappa_i,\qquad
		aI\preceq A_{\bm n}^{(i)}\preceq bI,
		\]
		where $0<a\le b<\infty$ are independent of $\bm n$. Then
		\begin{equation}
			\{G_{\bm n}\}_{\bm n}\sim_{\mathrm{GLT}}g,\quad \text{ where }\quad
			\begin{aligned}
				G_{\bm n}&=\Karcher(A_{\bm n}^{(1)},\ldots,A_{\bm n}^{(k)}),\\
				g(z)&=\Karcher(\kappa_1(z),\ldots,\kappa_k(z)).
			\end{aligned}
			\label{eq:glt-bounded-conclusion}
		\end{equation}
	\end{theorem}
	
	\begin{proof}
		By \ref{glt:distribution} and the spectral bounds, $aI_s\preceq\kappa_i(z)\preceq bI_s$ on a common full-measure subset $\Omega_0\subseteq\Omega$ \cite[Proposition~4.10, Theorem~4.11, and Corollary~2.17]{GLT20}.
		By \ref{km:continuity}, the pointwise mean $g$ is measurable on $\Omega_0$; define it arbitrarily on the zero-measure set.
		Choose
		\[
		X_{\bm n,0}=\frac1k\sum_iA_{\bm n}^{(i)},\qquad
		x_0=\frac1k\sum_i\kappa_i,
		\]
		and take $R,\hat{\theta},\rho$ as in Theorem~\ref{th:task}.
		Induction on $\nu$, using \ref{glt:algebra}--\ref{glt:functional-calculus} in \eqref{eq:R}, gives
		\[
		\{X_{\bm n,\nu}\}_{\bm n}\sim_{\mathrm{GLT}}x_\nu,
		\]
		where $x_\nu(z)$ obeys the same iteration with inputs $\kappa_i(z)$ and initial value $x_0(z)$. By Theorem~\ref{th:task}, all matrix and symbol iterates remain in $[aI,bI]$;
		for square roots, apply \ref{glt:functional-calculus} with the globally continuous function $t\mapsto\sqrt{\max\{t,0\}}$.
		The arguments of the logarithms have spectra in $[a/b,b/a]$, so the same rule applies to any globally continuous extension of the logarithm from that interval.
		We then apply \eqref{eq:frob} in dimension $m_{\bm n}$ and, for each admissible $z$, in dimension $s$. This gives
		\begin{equation}
			\frac{\norm{X_{\bm n,\nu}-G_{\bm n}}_2^2}{m_{\bm n}}
			\le b^2R^2\rho^\nu,\qquad
			\frac{\norm{x_\nu(z)-g(z)}_2^2}{s}
			\le b^2R^2\rho^\nu.
			\label{eq:glt-two-richardson-estimates}
		\end{equation}
		By \ref{acs:schatten}, the first inequality shows that $\{X_{\bm n,\nu}\}_{\bm n}$ is an a.c.s. for $\{G_{\bm n}\}_{\bm n}$.
		The second inequality gives
		\[
		\operatorname*{ess\,sup}_{z\in\Omega}\norm{x_\nu(z)-g(z)}
		\le\sqrt{s}\,bR\rho^{\nu/2}\longrightarrow0,
		\]
		and hence $x_\nu\to g$ in measure.
		The conclusion follows from \ref{glt:closure}.
	\end{proof}

	\subsection{Squeezing lemma for GLT sequences}
	To remove the uniform spectral bounds, we approximate the target from above and below in the L\"owner order. The order bounds also control its a.c.s. distance from the lower approximation: the singular-value gauge of this error is bounded by the gauge of the gap between the two approximants.
	Since those approximants are already GLT sequences, their gap can be estimated through their symbols. We datail such process in the following lemma.
	\begin{lemma}\label{lem:glt-squeeze}
		With $m_{\bm n}$ and $\Omega$ as in Section~\ref{sec:preliminaries}, let $L_{\bm n,r},H_{\bm n},U_{\bm n,r}$ be Hermitian matrices of size $m_{\bm n}$ such that
		\[
		L_{\bm n,r}\preceq H_{\bm n}\preceq U_{\bm n,r}
		\qquad\text{for every $\bm n$ and $r\in\N$}.
		\]
		Suppose that, for each fixed $r$,
		\[
		\{L_{\bm n,r}\}_{\bm n}\sim_{\mathrm{GLT}}\ell_r,\qquad
		\{U_{\bm n,r}\}_{\bm n}\sim_{\mathrm{GLT}}u_r,
		\]
		and that $\ell_r,u_r\to h$ in measure on $\Omega$ as $r\to\infty$,
		where $h:\Omega\to\C^{s\times s}$ is a matrix-valued measurable function.
		Then $\{H_{\bm n}\}_{\bm n}\sim_{\mathrm{GLT}}h$.
	\end{lemma}
	
	\begin{proof}
		Use $\psi(t)=t/(1+t)$ and $p_\psi$ from Section~\ref{sec:acs-prelim}.
		As stated in \eqref{eq:acs-pseudometric-prelim},
		\begin{equation}
			d_\psi(\{B_{\bm n}\},\{C_{\bm n}\})
			:=\limsup_{\bm n\to\infty}p_\psi(B_{\bm n}-C_{\bm n}),
			\label{eq:glt-acs-gauge}
		\end{equation}
		and as recalled there, convergence to zero in this pseudometric characterizes a.c.s. approximation. For GLT sequences with symbols $b,c$, \ref{glt:isometry} gives
		\begin{equation}
			d_\psi(\{B_{\bm n}\},\{C_{\bm n}\})
			=\frac1{s|\Omega|}\int_\Omega
			\sum_{\alpha=1}^s\psi(\sigma_\alpha(b(z)-c(z)))\,\dd z.
			\label{eq:glt-gauge-isometry}
		\end{equation}
		The bounded gauge identity in \ref{acs:gauge} applies even though $\psi$ is not compactly supported.
		
		For every fixed $r$, the order assumption gives
		\[
		O\preceq H_{\bm n}-L_{\bm n,r}
		\preceq U_{\bm n,r}-L_{\bm n,r}.
		\]
		Both differences are positive semidefinite, so their singular values coincides with the eigenvalues. The order comparison uses only the min--max principle and monotonicity of $\psi$, giving
		\[
		p_\psi(H_{\bm n}-L_{\bm n,r})
		\le p_\psi(U_{\bm n,r}-L_{\bm n,r}).
		\]
		Taking the limit superior in $\bm n$ and applying \ref{glt:isometry}, in the form \eqref{eq:glt-gauge-isometry}, only to the sequences $\{U_{\bm n,r}\}_{\bm n}$ and $\{L_{\bm n,r}\}_{\bm n}$, which are already known to be GLT, we obtain
		\begin{equation}
			\begin{aligned}
				d_\psi(\{H_{\bm n}\},\{L_{\bm n,r}\})
				&\le d_\psi(\{U_{\bm n,r}\},\{L_{\bm n,r}\})\\
				&=\frac1{s|\Omega|}\int_\Omega
				\sum_{\alpha=1}^s\psi(\sigma_\alpha(u_r(z)-\ell_r(z)))\,\dd z.
			\end{aligned}
			\label{eq:glt-squeeze-gauge}
		\end{equation}
		Since $D_r=u_r-\ell_r\to O_s$ in measure, the characterization in Section~\ref{sec:acs-prelim} gives
		\[
		p_\psi^{\mathrm{meas}}(D_r)
		=d_\psi^{\mathrm{meas}}(u_r,\ell_r)\longrightarrow0.
		\]
		Thus \eqref{eq:glt-squeeze-gauge} implies $d_\psi(\{H_{\bm n}\},\{L_{\bm n,r}\})\to0$.
		Consequently, $\{\{L_{\bm n,r}\}_{\bm n}\}_r$ is an a.c.s. for $\{H_{\bm n}\}_{\bm n}$. Since $\ell_r\to h$ in measure, \ref{glt:closure} gives $\{H_{\bm n}\}_{\bm n}\sim_{\mathrm{GLT}}h$.
	\end{proof}

	\subsection{GLT closure without uniform spectral bounds}
	We construct the two approximants required by the squeezing lemma by normalizing the inputs and their inverses separately. Adding a positive multiple of the identity places each normalized family in a common spectral interval, where the bounded closure theorem applies. Congruence invariance returns the means to the original scale, and self-duality
	provides the lower approximation.
	
	\begin{theorem}\label{thm:glt-general}
		Let $A_{\bm n}^{(i)}\succ O$ and
		$\{A_{\bm n}^{(i)}\}_{\bm n}\sim_{\mathrm{GLT}}\kappa_i$ for $i=1,\ldots,k$.
		Suppose that each matrix-valued symbol $\kappa_i$ is finite and HPD almost everywhere on $\Omega$. Then
		\begin{equation}
			\{G_{\bm n}\}_{\bm n}\sim_{\mathrm{GLT}}g,\qquad
			\begin{aligned}
				G_{\bm n}&=\Karcher(A_{\bm n}^{(1)},\ldots,A_{\bm n}^{(k)}),\\
				g&=\Karcher(\kappa_1,\ldots,\kappa_k).
			\end{aligned}
			\label{eq:glt-general-conclusion}
		\end{equation}
	\end{theorem}
	
	\begin{proof}
		Normalize both the inputs and their inverses. Define
		\begin{equation}
			\begin{aligned}
				S_{\bm n}&=\sum_i A_{\bm n}^{(i)},
				&\quad\sigma&=\sum_i\kappa_i,\\
				T_{\bm n}&=\sum_i(A_{\bm n}^{(i)})^{-1},
				&\quad\tau&=\sum_i\kappa_i^{-1},
			\end{aligned}
			\label{eq:glt-normalizers}
		\end{equation}
		so that $\{S_{\bm n}\}\sim_{\mathrm{GLT}}\sigma$ and $\{T_{\bm n}\}\sim_{\mathrm{GLT}}\tau$.
		Let 
		\begin{equation}
			{C_{\bm n}^{(i)}=S_{\bm n}^{-1/2}A_{\bm n}^{(i)}S_{\bm n}^{-1/2},
				\qquad
				D_{\bm n}^{(i)}=T_{\bm n}^{-1/2}(A_{\bm n}^{(i)})^{-1}T_{\bm n}^{-1/2},}
			\label{eq:glt-normalized-inputs}
		\end{equation}
		where both terms are well-defined; moreover, the corresponding GLT sequences have respective symbols
		\[
		c_i=\sigma^{-1/2}\kappa_i\sigma^{-1/2},\qquad
		d_i=\tau^{-1/2}\kappa_i^{-1}\tau^{-1/2}.
		\]
		Since $A_{\bm n}^{(i)}\preceq S_{\bm n}$ and $(A_{\bm n}^{(i)})^{-1}\preceq T_{\bm n}$, the preservation of order under congruence gives
		\[
		O\prec C_{\bm n}^{(i)}\preceq I,\qquad
		O\prec D_{\bm n}^{(i)}\preceq I.
		\]
		If $\varepsilon>0$ is fixed, both sets of matrices $C_{\bm n}^{(i)}+\varepsilon I$ and $D_{\bm n}^{(i)}+\varepsilon I$ therefore lie in the uniform interval $[\varepsilon I,(1+\varepsilon)I]$.
		By Theorem~\ref{thm:glt-bounded}, the following Karcher mean sequences are GLT sequences:
		\[
		\begin{aligned}
			P_{\bm n,\varepsilon}
			&=\Karcher(C_{\bm n}^{(1)}+\varepsilon I,\ldots,
			C_{\bm n}^{(k)}+\varepsilon I),\\
			Q_{\bm n,\varepsilon}
			&=\Karcher(D_{\bm n}^{(1)}+\varepsilon I,\ldots,
			D_{\bm n}^{(k)}+\varepsilon I),
		\end{aligned}
		\]
		with symbols 
		\[
		\begin{aligned}
			p_\varepsilon&=\Karcher(c_1+\varepsilon I_s,\ldots,c_k+\varepsilon I_s),\\
			q_\varepsilon&=\Karcher(d_1+\varepsilon I_s,\ldots,d_k+\varepsilon I_s).
		\end{aligned}
		\]
		By \ref{km:monotonicity} and \ref{km:commuting} we infer
		\[
		\varepsilon I\preceq P_{\bm n,\varepsilon},Q_{\bm n,\varepsilon}
		\preceq(1+\varepsilon)I.
		\]
		We next use congruence invariance to return to the original scale.
		Let
		\begin{equation}
			U_{\bm n,\varepsilon}
			=S_{\bm n}^{1/2}P_{\bm n,\varepsilon}S_{\bm n}^{1/2},
			\qquad
			L_{\bm n,\varepsilon}
			=T_{\bm n}^{-1/2}Q_{\bm n,\varepsilon}^{-1}T_{\bm n}^{-1/2}.
			\label{eq:glt-bounds-definition}
		\end{equation}
		so that the corresponding sequences are GLT. We apply
		\ref{km:congruence} with $M=S_{\bm n}^{1/2}$ and
		$M=T_{\bm n}^{1/2}$ obtaining
		\begin{align}
			U_{\bm n,\varepsilon}
			&=\Karcher(A_{\bm n}^{(1)}+\varepsilon S_{\bm n},\ldots,
			A_{\bm n}^{(k)}+\varepsilon S_{\bm n}),
			\label{eq:glt-upper-congruence}\\
			T_{\bm n}^{1/2}Q_{\bm n,\varepsilon}T_{\bm n}^{1/2}
			&=\Karcher((A_{\bm n}^{(1)})^{-1}+\varepsilon T_{\bm n},\ldots,
			(A_{\bm n}^{(k)})^{-1}+\varepsilon T_{\bm n}).
			\label{eq:glt-lower-congruence}
		\end{align}
		Inverting \eqref{eq:glt-lower-congruence} gives
		\begin{equation}
			L_{\bm n,\varepsilon}
			=\Karcher((A_{\bm n}^{(1)})^{-1}+\varepsilon T_{\bm n},\ldots,
			(A_{\bm n}^{(k)})^{-1}+\varepsilon T_{\bm n})^{-1}.
			\label{eq:glt-lower-identity}
		\end{equation}
		Now, by \ref{km:monotonicity} and \eqref{eq:glt-upper-congruence},
		$G_{\bm n}\preceq U_{\bm n,\varepsilon}$; for the lower bound, \ref{km:monotonicity} and \ref{km:self-duality} give
		\[
		\begin{aligned}
			&\Karcher((A_{\bm n}^{(1)})^{-1}+\varepsilon T_{\bm n},\ldots,
			(A_{\bm n}^{(k)})^{-1}+\varepsilon T_{\bm n})\\
			&\hspace{2em}\succeq
			\Karcher((A_{\bm n}^{(1)})^{-1},\ldots,(A_{\bm n}^{(k)})^{-1})
			=G_{\bm n}^{-1}.
		\end{aligned}
		\]
		Inversion reverses the order, so \eqref{eq:glt-lower-identity} gives
		\begin{equation}
			L_{\bm n,\varepsilon}\preceq G_{\bm n}\preceq U_{\bm n,\varepsilon}.
			\label{eq:glt-loewner-sandwich}
		\end{equation}
		By \ref{glt:algebra}--\ref{glt:functional-calculus}, the symbols in
		\eqref{eq:glt-bounds-definition} are
		\[
		u_\varepsilon=\sigma^{1/2}p_\varepsilon\sigma^{1/2},\qquad
		\ell_\varepsilon=\tau^{-1/2}q_\varepsilon^{-1}\tau^{-1/2},
		\]
		and applying \ref{km:congruence} in dimension $s$ gives
		\begin{equation}
			\begin{aligned}
				u_\varepsilon
				&=\Karcher(\kappa_1+\varepsilon\sigma,\ldots,
				\kappa_k+\varepsilon\sigma),\\
				\ell_\varepsilon
				&=\Karcher(\kappa_1^{-1}+\varepsilon\tau,\ldots,
				\kappa_k^{-1}+\varepsilon\tau)^{-1}.
			\end{aligned}
			\label{eq:glt-bound-symbols}
		\end{equation}
		For almost every $z$, the inputs and their inverses are finite positive definite matrices. By \ref{km:continuity}, continuity of inversion and \ref{km:self-duality},
		\begin{equation}
			u_\varepsilon(z)\longrightarrow g(z),\qquad
			\ell_\varepsilon(z)\longrightarrow
			\Karcher(\kappa_1(z)^{-1},\ldots,\kappa_k(z)^{-1})^{-1}
			=g(z).
			\label{eq:glt-symbol-squeeze-limit}
		\end{equation}
		Since $|\Omega|<\infty$, this convergence also holds in measure. Apply Lemma~\ref{lem:glt-squeeze} to \eqref{eq:glt-loewner-sandwich}, with $\varepsilon=1/r$. This gives \eqref{eq:glt-general-conclusion} and proves the theorem.
	\end{proof}
	
	\noindent The spectral distribution follows from \ref{glt:distribution} because each $G_{\bm n}$ is Hermitian; for every continuous compactly supported function $F:\mathbb R\to\C$, we have
	\begin{equation}
		\lim_{\bm n\to\infty}\frac1{sN(\bm n)}
		\sum_{j=1}^{sN(\bm n)}F(\lambda_j(G_{\bm n}))\\
		=
		\frac1{s|\Omega|}\int_\Omega
		\sum_{\alpha=1}^sF(\lambda_\alpha(g(z)))\,\dd z.
		\label{eq:glt-karcher-spectral-distribution}
	\end{equation}
	By \ref{km:commuting}, if the symbols commute pairwise almost everywhere,
	\[
	g=\exp\!\left(\frac1k\sum_i\log\kappa_i\right)
	=(\kappa_1\cdots\kappa_k)^{1/k},
	\]
	and in particular, this formula applies when $s=1$.
	
	{We remark here that the} $\varepsilon$-approximation idea can be traced back to \cite{DEG25}, where the maximal result for two-input GLT matrix sequences was obtained improving the analysis in \cite{GEO25} and giving a complete proof of Conjecture 10.1 in \cite{GS17}.
	
	{Moreover, the $s$-block $d$-level GLT matrix sequence have the rigidity of the prescribed matrix size $m_{\bm n}=sN(\bm n)$. This is reasonable when approximating partial differential equations on Cartesian domains: however, even in that case, the presence of boundary conditions or of high order methods can change the size to $s(N(\bm n)+o(1))$. We observe that the latter is not a limitation due to the combination of the extra-dimensional approach and of the reduced versions of the GLT theory so that any increasing matrix size $d_{\bm n}$ is allowed, by maintaining all the machinery of the standard GLT theory; see \cite{reduced,extra-dim} and references therein.

		\section{Experiments on spectral distributions}\label{sec:numerical-glt}
		In this section we propose a numerical assessment of the results examined in the current paper.
		In is worth noticing that by definition \eqref{eq:karcher-variational} the Karcher mean of $k$ HPD matrices is HPD itself, so that, when considering the associated matrix-sequence, its GLT symbol naturally translates to the eigenvalue distribution.
		
		In the following numerical experiments, we test the {quality} of the approximation provided by the eigenvalue symbol, under different parameters.

		\paragraph{{Univariate scalar symbol.}}
		We begin numerical considerations of the theoretical results discussed in a unilevel scalar setting. In particular, we examine the shrinking ratio to the theoretical spectral distribution 
		\begin{itemize}
			\item of the approximation provided by \eqref{eq:R} for increasing number of iteration while keeping the dimension fixed;
			\item of the matrix-sequence given by the Karcher mean for increasing size and a fixed number of iterations.
		\end{itemize}
		Consider the following finite and a.e. positive functions {
			\begin{displaymath}
				a(x) = \left\{
				\begin{array}{cl}
					x^2 + 2 &\qquad \text{for } x\in\left[0, \frac{1}{2}\right]\\[5pt]
					\log\left(x+ \frac{1}{2}\right) &\qquad \text{for } x\in\left(\frac{1}{2},1\right]
				\end{array}\right. ,
				\qquad\quad
				\begin{aligned}
					f_1(\theta) &= 9 - 4\cos(\theta) - 4\cos(2\theta),\\
					f_2(\theta) &= \theta^2,
				\end{aligned}
			\end{displaymath}
			with $\theta\in\left[ -\pi, \pi \right]$. Let} $D_n(a)$, $T_n(f_1)$ and $T_n(f_2)$ be the associated diagonal sampling and Toeplitz matrices, and set $G_n = \mathcal{G}\left(D_n(a), T_n(f_1), T_n(f_2)\right)$, $g\left(x,\theta\right) = \mathcal{G}\left(a(x), f_1(\theta), f_2(\theta)\right)$. 
		As a consequence of Theorem \ref{thm:glt-general}, we deduce
		\begin{displaymath}
			\left\{G_n\right\}_n \sim_{ \lambda } g \left( x, \theta\right)\quad
			\text{in } \left[ 0, 1 \right] \times \left[ -\pi, \pi \right];
		\end{displaymath}
		since $g$ is a.e. continuous on its domain, we have that the spectrum of $G_n$ provides an approximate reconstruction of the symbol for $n$ large enough.
		
		In Figure \ref{fig:esp_1} {we show this behavior}. {In the plot on the left we fix the amount of Richardson iterations as $\nu = 10$ and check the adherence of the eigenvalues of $X_{n,\nu}$ to the symbol $g$, while in the plot on the right} we consider a fixed size $m_n=200$ and {show} the  {eigenvalues} of $X_{n,\nu}$ for $\nu = 2, 4, 10$.
		{In these tests we fixed both the Richardson step size to $\hat{\theta}=0.11$ and the initial guess to $X_{n,0} = 3I_n$.
			As stated in} the aforementioned theoretical results, we observe that said spectrum shrinks on the graph of the symbol $g$ as $\nu$ grows.
		
		\begin{figure}
			\centering
			\includegraphics[width=1.0\textwidth]{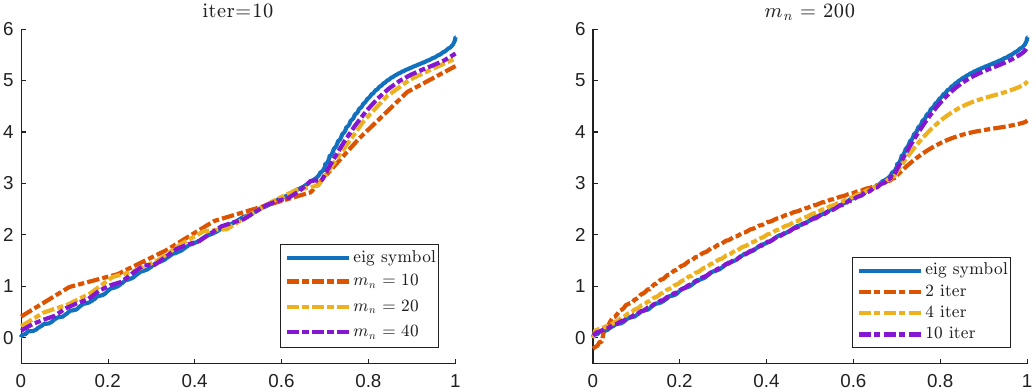}
			\caption{Convergence to the theoretical distribution (full line) of the spectra (dashed lines) of Karcher means for increasing size $m_n$ {(left)} and Richardson iteration number {(right)} for the univariate scalar test.}
			\label{fig:esp_1} 
		\end{figure}
		
		\paragraph{The choice of step size and initial guess in the scalar case}
		{We perform a further experiment in the univariate scalar case:} we analyze how the choice of the step size and initial guess influences the convergence speed of the Richardson algorithm.
		We employ the same matrices of the previous example and fix the size {to} $m_n=200$.
		
		{Firstly}, we fixed the initial guess $X_{n,0} = 3I_n$ and consider different step sizes.
		Heuristic observations lead to the conclusion that the iteration \eqref{eq:R} is extremely sensible to the step size, and, when compared to a potentially large number of iterations, a small $\hat{\theta}$ is heavily preferable {in this specific setting}.
		In this test, we set $\hat{\theta} = 0.07, 0.10, 0.13$ and {the plot on the left of Figure \ref{fig:esp_2} shows that} a larger $\hat{\theta}$ fails to achieve acceptable convergence {near 0} even with a sensibly more substantial number of iterations.
		
		{On the other hand, we fixed step size $\hat{\theta}=0.10$ and consider the following initial guesses
			\begin{displaymath}
				X_{n,0}^{\left(1\right)} = 3I_n,\quad
				X_{n,0}^{\left(2\right)} = M_n,\quad
				X_{n,0}^{\left(3\right)} = H_n,\quad
				X_{n,0}^{\left(4\right)} = K_n,
		\end{displaymath}}
		where $M_n$ is the arithmetic mean of the matrices involved and 
		\begin{displaymath}
			{H_n = D_n(a^{1/6}) T_n(f_1)^{1/3} T_n(f_2)^{1/3} D_n(a^{1/6}),
				\qquad
				K_n = D_n(a^{1/6}) C_n(f_1)^{1/3} C_n(f_2)^{1/3} D_n(a^{1/6}).}
		\end{displaymath}
		Here $C_n(f_1)$, $C_n(f_2)$ are the circulant matrices with spectra given by {samplings} of $f_1$ and $f_2$ respectively on a uniform grid over $\left[ -\pi, \pi \right]$. {Since $f_2(0) = 0$, its evaluation at the specific circulant grid node $0$ is replaced by its value on the next node, thus avoiding issues regarding the singularity of $C_n(f_2)$: this trick is classical and leads to a Strang-type circulant approximation \cite{CN}.}
		With such choice, $\left\{H_n\right\}_n, \left\{K_n\right\}_n \sim_{ \mathrm{GLT}} g $ in $ \left[0, 1\right] \times \left[-\pi, \pi\right]$, i.e. the initialization guess for the Richardson algorithm shares the same GLT distribution of its limit point $\left\{G_n\right\}_n$.
		It is worth noticing that the formula for $H_n$ requires the explicit computation of the SVD of $T_n(f_1)$ and $T_n(f_2)$, while $K_n$ can be assembled with $\mathcal{O}(m_n\log_2 m_n)$ operations.
		The resulting spectral distribution is displayed in the {plot on the right of Figure \ref{fig:esp_2}}.
		
		We recognize that the Toeplitz GLT guess $ H_n$ achieves much better performances than its circulant counterpart $ K_n$, while the other choices appear to be roughly equivalent.
		
		\begin{figure}
			\centering
			\includegraphics[width=\textwidth]{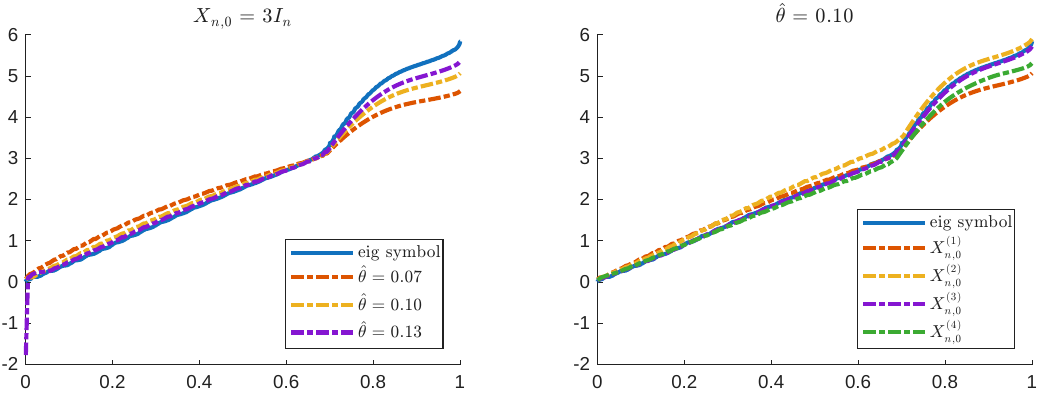}
			\caption{Convergence to the theoretical distribution (full line) of the spectra (dashed lines) of Karcher means for different {step size (left) and different initial guesses (right)} for the univariate scalar test. Here $m_n = 200$, while the iteration number is set at 5.}
			\label{fig:esp_2} 
		\end{figure}

		\paragraph{Matrix-valued symbols.}
		For our third experiment, we switch the framework from scalar to matrix-valued symbols. Consider the following a.e. HPD matrix-valued functions
		{\begin{displaymath}
				a(x) = \begin{bmatrix}
					2 & e^{\imagy\pi x}\\
					e^{-\imagy\pi x} &2
				\end{bmatrix},
				\quad
				F_1(\theta) = \begin{bmatrix}
					3 & -\cos(\theta)\\
					-\cos(\theta) &3
				\end{bmatrix},
				\quad
				F_2(\theta) = \begin{bmatrix}
					1 & 1-\frac{\abs{\theta}}{\pi} \\
					1 - \frac{\abs{\theta}}{\pi} & 1
				\end{bmatrix},
			\end{displaymath}
			where $x\in\left[0, 1\right]$ and $\theta\in\left[ -\pi, \pi \right]$.}
		
		Let $D_n(a)$, $T_n(F_1)$ and $T_n(F_2)$ the associated block diagonal sampling matrix and block Toeplitz matrices.
		Note that both $D_n(a)$ and $T_n(F_1)$ are well-conditioned, while $T_n(F_2)$ has its minimal eigenvalue approaching zero as $1/m_n$.
		Consider $G_n=\mathcal{G}\left(D_n(a), T_n(F_1), T_n(F_2)\right)$, and set $g(x,\theta) = \mathcal{G}\left(a(x), F_1(\theta), F_2(\theta)\right)$. \{Again, by} Theorem \ref{thm:glt-general} and definition \eqref{eq:karcher-variational}, $\left\{ G_n\right\}_n \sim_\lambda g$ in $\left[0, 1\right]\times\left[-\pi, \pi\right]$.
	
\begin{figure}
		\centering
		\includegraphics[width=1.0\textwidth]{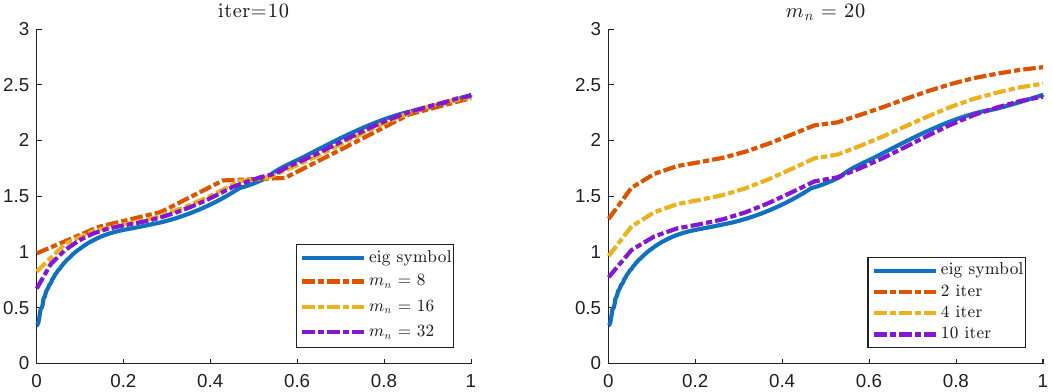}
		\caption{Convergence to the theoretical distribution (full line) of the spectra (dashed lines) of Karcher means for increasing size $m_n$ {(left)} and Richardson iteration number {(right) for the matrix-valued symbol test.}
			The initial guess is set at $3 I_n$ and $\hat{\theta} = 0.10$.}
		\label{fig:esp_31} 
	\end{figure}

	{In Figure \ref{fig:esp_31} we report the results of the analogous test done for the univariate scalar case, that is, the adherence of the eigenvalues of $G_n$ to the spectral symbol $g$;} moreover, as inferred from the {plot on the right}, the fast-converge properties of the Richardson algorithm are preserved in this setting.
	
\begin{figure}
		\centering
		\includegraphics[width=1.0\textwidth]{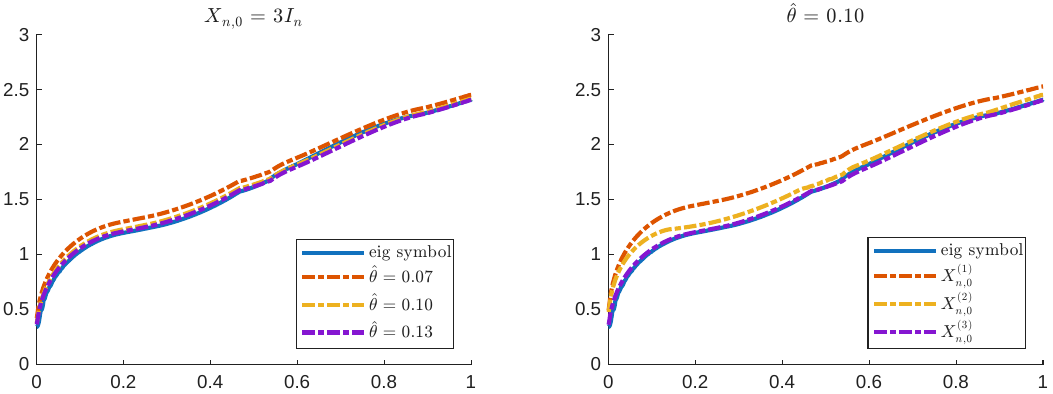}
		\caption{Convergence to the theoretical distribution (full line) of the spectra (dashed lines) of Karcher means for different {step sizes (left)} and {initial guesses (right) for the matrix-valued symbol test}. The matrices dimension is fixed at $m_n = 200$ and {the Richardson iterations counter at $10$ for the plot on the left, while at $4$ for the plot on the right.}}
		\label{fig:esp_32} 
	\end{figure}

	{Moreover}, we check whether the choice of the initial guess and sensibility {with respect} to the step size are still valid for matrices with matrix-valued symbols; results in Figure \ref{fig:esp_32} are in strong analogy with what presented before and the method is again extremely sensible to the step size, {while the GLT initialization provides better performances even for a very low amount of iterations}.	

\begin{figure}
		\centering
		\includegraphics[width=1.0\textwidth]{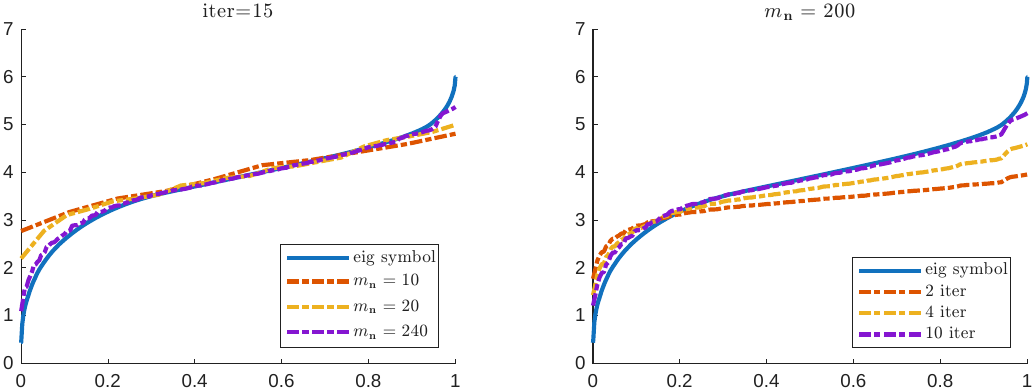}
		\caption{{Convergence to the theoretical distribution (full line) of the spectra (dashed lines) of Karcher means for increasing size $m_{\boldsymbol{n}}$ (left) and Richardson iteration number (right) for the multivariate test.}}
		\label{fig:esp_41}
	\end{figure}

	\paragraph{Multivariate symbols.}
	{Our} final experiment {explores the} multi-level setting {(now $n$ is to be considered as a multi-index $\boldsymbol{n}$).}
	Consider the following functions {for $x\in\left[ 0, 1 \right]$ and $\theta_1, \theta_2 \in \left[ -\pi, \pi \right]$}
	\begin{displaymath}
		a(x) = e^x,\qquad
		f_1(\theta_1, \theta_2) = 4 - 2\cos(\theta_1) - 2\cos(\theta_2),\qquad
		f_2(\theta_1, \theta_2) = 10 + \theta_1^2 \cos(\theta_2),
	\end{displaymath}
	and the associated diagonal sampling and Toeplitz matrices; note that, since its minimal eigenvalue approaches zero from above as $1/m_{\boldsymbol{n}}^2$, $T_{\boldsymbol{n}}(f_1)$ is ill-conditioned.
	In Figure \ref{fig:esp_41} {we show the adherence} of the {eigenvalues of} the geometric mean $G_n$ {to its} theoretical distribution {as done in the previous tests}.
	The initialization matrix is set in both cases as $3I_{\boldsymbol{n}}$, well inside the theoretical spectral bounds.
	
\begin{figure}
		\centering
		\includegraphics[width=1.0\textwidth]{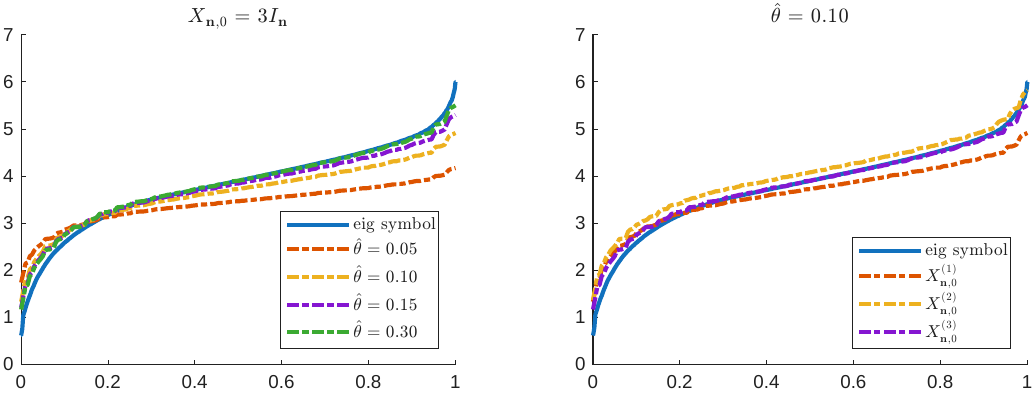}
		\caption{{Convergence to the theoretical distribution (full line) of the spectra (dashed lines) of Karcher means for different step sizes (left) and initial guesses (right) for the multivariate test.}}
		\label{fig:esp_42} 
	\end{figure}

	{Moreover, in the same fashion as before}, in Figure \ref{fig:esp_42} we show the results for different step sizes and initialization points.
	As before, $X_{\boldsymbol{n},0}^{(1)} = 3I_{\boldsymbol{n}}$, $X_{\boldsymbol{n},0}^{(2)}$ {is} the arithmetic mean and
	\begin{displaymath}
		X_{\boldsymbol{n},0}^{(3)} = 
		D_{\boldsymbol{n}}(a^{1/6}) T_{\boldsymbol{n}}(f_1)^{1/3}
		T_{\boldsymbol{n}}(f_2)^{1/3} D_{\boldsymbol{n}}(a^{1/6}).
	\end{displaymath}
	{Similar} considerations {outlined in} the previous sections apply here {as well}. 
	However, it is worth noticing that, in contrast with the other examples, is this case a larger step size appears to be able to accelerate the convergence.
	

	\section{Conclusions}\label{sec:end}
	In the current work, we have obtained normalized Schatten estimates independent of the matrix dimension, in the context of the additive Richardson iteration for the computation of the Karcher mean. These estimates have been employed for proving that Karcher means of structured matrices are asymptotically structured, both in the Toeplitz and GLT settings, and with weak assumption on the input GLT symbols. 
	Ad hoc numerical experiments have been reported and discussed. As future lines of research we can mention the following items.
	\begin{itemize}
		\item The GLT machinery is quite general and hence other iterative methods for computing alternative geometric means can be considered.
		\item It would be interesting to provide quantitative estimates on the computation gain, when choosing as starting point a GLT matrix sequence having as GLT symbol already the exact geometric mean of the imput symbols.
		\item A more systematic numerical investigation is required for general $r,d>1$, with special attention to cases stemming from concrete applications.
	\end{itemize}
	
	\paragraph{Acknowledgements.} The conjecture in \cite[pp. 13-14]{GEO25} was firstly discussed ten years before by the third author with Giovanni Barbarino, Carlo Garoni, Bruno Iannazzo, Debora Sesana.


\begin{thebibliography}{99}
		\urlstyle{same}
		\setlength{\itemsep}{2pt}
		
		\bibitem{extra-dim}
		A. Adriani, A. J. A. Schiavoni-Piazza, S. Serra-Capizzano, and C. Tablino-Possio, \emph{Revisiting the notion of approximating class of sequences for handling approximated PDEs on moving or unbounded domains}, Electronic Transactions on Numerical Analysis \textbf{63} (2025), 424--451.
		\href{https://doi.org/10.1553/etna_vol63s424}{doi:\nolinkurl{10.1553/etna_vol63s424}}.
		
		\bibitem{GEO25}
		D. Ahmad, M. F. Khan, and S. Serra-Capizzano, \emph{Matrix-sequences of geometric means in the case of hidden (asymptotic) structures}, Mathematics \textbf{13} (2025), no. 3, article 393.
		\href{https://doi.org/10.3390/math13030393}{doi:\nolinkurl{10.3390/math13030393}}.
		
		\bibitem{ALM04}
		T. Ando, C.-K. Li, and R. Mathias, \emph{Geometric means}, Linear Algebra and its Applications \textbf{385} (2004), 305--334. \href{https://doi.org/10.1016/j.laa.2003.11.019}{doi:\nolinkurl{10.1016/j.laa.2003.11.019}}.
		
		\bibitem{reduced}
		G. Barbarino, \emph{A systematic approach to reduced GLT}, BIT \textbf{62-3} (2022), 681--743.
		\href{https://doi.org/10.1007/s10543-021-00896-7}{doi:\nolinkurl{10.1007/s10543-021-00896-7}}.
		
		\bibitem{BG17}
		G. Barbarino and C. Garoni, \emph{From convergence in measure to convergence of matrix-sequences through concave functions and singular values}, Electronic Journal of Linear Algebra \textbf{32} (2017), 500--513.
		\href{https://doi.org/10.13001/1081-3810.3663}{doi:\nolinkurl{10.13001/1081-3810.3663}}.
		
		\bibitem{GLT20-uni}
		G. Barbarino, C. Garoni, and S. Serra-Capizzano, \emph{Block generalized locally Toeplitz sequences: theory and applications in the unidimensional case}, Electronic Transactions on Numerical Analysis \textbf{53} (2020), 28--112.
		\href{https://doi.org/10.1553/etna_vol53s28}{doi:\nolinkurl{10.1553/etna_vol53s28}}.
		
		\bibitem{GLT20}
		G. Barbarino, C. Garoni, and S. Serra-Capizzano, \emph{Block generalized locally Toeplitz sequences: theory and applications in the multidimensional case}, Electronic Transactions on Numerical Analysis \textbf{53} (2020), 113--216.
		\href{https://doi.org/10.1553/etna_vol53s113}{doi:\nolinkurl{10.1553/etna_vol53s113}}.
		
		\bibitem{BH06}
		R. Bhatia and J. Holbrook, \emph{Riemannian geometry and matrix geometric means}, Linear Algebra and its Applications \textbf{413} (2006), nos.~2--3, 594--618.
		\href{https://doi.org/10.1016/j.laa.2005.08.025}{doi:\nolinkurl{10.1016/j.laa.2005.08.025}}.
		
		\bibitem{BI13}
		D. A. Bini and B. Iannazzo, \emph{Computing the Karcher mean of symmetric positive definite matrices}, Linear Algebra and its Applications \textbf{438} (2013), no.~4, 1700--1710.
		\href{https://doi.org/10.1016/j.laa.2011.08.052}{doi:\nolinkurl{10.1016/j.laa.2011.08.052}}.
		
		\bibitem{BGr}
		A. B\"ottcher and S. Grudsky, \emph{Spectral properties of banded Toeplitz matrices}, Society for Industrial and Applied Mathematics (SIAM), Philadelphia, PA, 2005.
		\href{https://doi.org/10.1137/1.9780898717853}{doi:\nolinkurl{10.1137/1.9780898717853}}.
		
		\bibitem{CN}
		R. Chan and M. K. Ng, \emph{Conjugate gradient methods for Toeplitz systems}, SIAM Review \textbf{38-3} (1996), 427--482.
		\href{https://doi.org/10.1137/S0036144594276474}{doi:\nolinkurl{10.1137/S0036144594276474}}.
		
		\bibitem{GS17}
		C. Garoni and S. Serra-Capizzano, \emph{Generalized locally Toeplitz sequences: theory and applications}, vol. I, Springer, Cham, 2017.
		\href{https://doi.org/10.1007/978-3-319-53679-8}{doi:\nolinkurl{10.1007/978-3-319-53679-8}}.
		
		\bibitem{GS18}
		C. Garoni and S. Serra-Capizzano, \emph{Generalized locally Toeplitz sequences: theory and applications}, vol.~II, Springer, Cham, 2018.
		\href{https://doi.org/10.1007/978-3-030-02233-4}{doi:\nolinkurl{10.1007/978-3-030-02233-4}}.
		
		\bibitem{IP18}
		B. Iannazzo and M. Porcelli, \emph{The Riemannian Barzilai--Borwein method with nonmonotone line search and the matrix geometric mean computation}, IMA Journal of Numerical Analysis \textbf{38} (2018), no.~1, 495--517.
		\href{https://doi.org/10.1093/imanum/drx015}{doi:\nolinkurl{10.1093/imanum/drx015}}.
		
		\bibitem{DEG25}
		A. Ilyas, M. F. Khan, V. Loi, and S. Serra-Capizzano, \emph{A maximal distribution result and numerical tests for geometric means of HPD GLT matrix-sequences with degenerate commuting/non-commuting GLT symbols}, Linear Algebra and its Applications (2025).
		\href{https://doi.org/10.1016/j.laa.2025.08.004}{doi:\nolinkurl{10.1016/j.laa.2025.08.004}}.
		
		\bibitem{TH}
		M. F. Khan, \emph{GLT matrix-sequences and few emblematic applications}, Ph.D. thesis, University of Insubria, 2025.
		\url{https://hdl.handle.net/11383/2201451}.
		
		
		\bibitem{LL11}
		J. Lawson and Y. Lim, \emph{Monotonic properties of the least squares mean}, Mathematische Annalen \textbf{351} (2011), no.~2, 267--279.
		\href{https://doi.org/10.1007/s00208-010-0603-6}{doi:\nolinkurl{10.1007/s00208-010-0603-6}}.
		
		\bibitem{LL14}
		J. Lawson and Y. Lim, \emph{Karcher means and Karcher equations of positive definite operators}, Transactions of the American Mathematical Society, Series B \textbf{1} (2014), no.~1, 1--22.
		\href{https://doi.org/10.1090/S2330-0000-2014-00003-4}{doi:\nolinkurl{10.1090/S2330-0000-2014-00003-4}}.
		
		\bibitem{LP12}
		Y. Lim and M. P\'alfia, \emph{Matrix power means and the Karcher mean}, Journal of Functional Analysis \textbf{262} (2012), no.~4, 1498--1514.
		\href{https://doi.org/10.1016/j.jfa.2011.11.012}{doi:\nolinkurl{10.1016/j.jfa.2011.11.012}}.
		
		\bibitem{LP18}
		Y. Lim and M. P\'alfia, \emph{Existence, uniqueness and an ODE approach to the $L^1$-Karcher mean}, Advances in Mathematics \textbf{376} (2021), article~107435.
		\href{https://doi.org/10.1016/j.aim.2020.107435}{doi:\nolinkurl{10.1016/j.aim.2020.107435}}.
		
		\bibitem{Moakher05}
		M. Moakher, \emph{A differential geometric approach to the geometric mean of symmetric
			positive-definite matrices}, SIAM Journal on Matrix Analysis and Applications \textbf{26} (2005), no.~3, 735--747.
		\href{https://doi.org/10.1137/S0895479803436937}{doi:\nolinkurl{10.1137/S0895479803436937}}.
		
		\bibitem{NgBook} M. K. Ng, \emph{Iterative methods for Toeplitz systems}, Oxford University Press, New York, 2004.
		\href{https://doi.org/10.1093/oso/9780198504207.003.0002}{doi:\nolinkurl{10.1093/oso/9780198504207.003.0002}}.
		
		\bibitem{KoroSINUM} S. Serra-Capizzano, \emph{A Korovkin-based approximation of multilevel Toeplitz matrices (with rectangular unstructured blocks) via multilevel trigonometric matrix spaces}, SIAM Journal on Numerical Analysis \textbf{36-6} (1999), 1831--1857.
		\href{https://doi.org/10.1137/S0036142997322497}{doi:\nolinkurl{10.1137/S0036142997322497}}.
		
		\bibitem{Sbit} S. Serra-Capizzano, \emph{Spectral and computational analysis of block Toeplitz matrices having nonnegative definite matrix-valued generating functions}, BIT \textbf{39-1} (1999), 152--175.
		\href{https://doi.org/10.1023/A:1022329526925}{doi:\nolinkurl{10.1023/A:1022329526925}}.
		
		\bibitem{acs} S. Serra-Capizzano, \emph{Distribution results on the algebra generated by Toeplitz sequences: a finite-dimensional approach}, Linear Algebra and its Applications \textbf{238-1/3} (2001), 121--130.
		\href{https://doi.org/10.1016/S0024-3795(00)00311-6}{doi:\nolinkurl{10.1016/S0024-3795(00)00311-6}}.
		
		\bibitem{GLT03} S. Serra Capizzano, \emph{Generalized locally Toeplitz sequences: spectral analysis and applications to discretized partial differential equations}, Linear Algebra and its Applications \textbf{366} (2003), 371--402.
		\href{https://doi.org/10.1016/S0024-3795(02)00504-9}{doi:\nolinkurl{10.1016/S0024-3795(02)00504-9}}.
		
		\bibitem{GLT06} S. Serra Capizzano, \emph{The GLT class as a generalized Fourier analysis and applications}, Linear Algebra and its Applications \textbf{419-1} (2006), 180--233.
		\href{https://doi.org/10.1016/j.laa.2006.04.012}{doi:\nolinkurl{10.1016/j.laa.2006.04.012}}.
		
		\bibitem{Tilli98}
		P.~Tilli, \emph{Locally Toeplitz sequences: spectral properties and applications}, Linear Algebra and its Applications \textbf{278} (1998), nos.~1--3, 91--120.
		\href{https://doi.org/10.1016/S0024-3795(97)10079-9}{doi:\nolinkurl{10.1016/S0024-3795(97)10079-9}}.
		
		\bibitem{Ty} E. E. Tyrtyshnikov, \emph{A unifying approach to some old and new theorems on distribution and clustering}, Linear Algebra and its Applications \textbf{232} (1996), 1--43.
		\href{https://doi.org/10.1016/0024-3795(94)00025-5}{doi:\nolinkurl{10.1016/0024-3795(94)00025-5}}.
		
	\end{thebibliography}
\end{document}